\documentclass[a4paper,10pt]{article}
\usepackage{amsmath,amssymb,amsthm,dsfont}
\usepackage[english]{babel}
\newtheorem{theo}{Theorem}[section] 
\newtheorem{defi}[theo]{Definition}
\newtheorem{lemm}[theo]{Lemma} 
\newtheorem{prop}[theo]{Proposition}
\newtheorem{coro}[theo]{Corollary}

\newcommand{\Na}{\mathbb N}                   

\newcommand{\Ra}{\mathbb R}                   

\newcommand{\Ca}{\mathbb C}                   

\newcommand{\scal}[1]{\langle #1 \rangle}

\newcommand{\finpreuve}{\hfill $\Box$}

\newcommand{\name}{$\underline{\qquad \qquad}$} \newcommand{\etan}{$ \& $ }

\begin{document}

\author{  Jean-Marc Bouclet \etan Julien Royer \\
Institut de Math\'ematiques de Toulouse \\ 
118 route de Narbonne 
\\ 
F-31062 Toulouse Cedex 9 \\
 jean-marc.bouclet@math.univ-toulouse.fr \\
 julien.royer@math.univ-toulouse.fr}
\title{Sharp low frequency resolvent estimates on asymptotically conical manifolds}

\date{ } 

\maketitle

\begin{abstract} On  a class of asymptotically conical manifolds, we prove two types of low frequency estimates for the resolvent of the Laplace-Beltrami operator. The first result is a uniform $ L^2 \rightarrow L^2 $ bound for $ \scal{r}^{-1} (- \Delta_G - z)^{-1} \scal{r}^{-1}  $ when $ \mbox{Re}(z) $ is small, with the optimal weight $ \scal{r}^{-1} $.  The second one is about powers of the resolvent. For any integer $N$, we prove uniform $ L^2 \rightarrow L^2 $ bounds for $ \scal{\epsilon r}^{-N} (-\epsilon^{-2} \Delta_G - Z)^{-N} \scal{\epsilon r}^{-N}  $ when $ \mbox{Re}(Z) $ belongs to a compact subset of $ (0,+\infty) $ and $ 0 < \epsilon \ll 1 $. These results are obtained by proving similar estimates on a pure cone with a long range perturbation of the metric at infinity.
\end{abstract}



\section{Introduction and main results}

The long range scattering theory of the Laplace-Beltrami operator on asymptotically Euclidean or conical manifolds has been widely studied. It has reached a point where our global understanding of the spectrum, in particular the behaviour of the resolvent at low, medium and high frequencies, allows to extend to curved settings many results which are well known on $ \Ra^n $.   We have typically in mind global in time Strichartz estimates \cite{Tataru1,MaMeTaru,MeTaru,HaZh,Zhang} or various instances of the local energy decay \cite{BonyHafner0,BonyHafner2,VasyWunschHardy,VasyWunsch2,BoucletCPDE,RoTa,BoucletRoyer} which are important tools in nonlinear PDE arising in mathematical physics. 
We refer to the recent paper \cite{RoTa} which surveys resolvent estimates (or limiting absorption principle) and some of their applications in this geometric framework.

In this picture, the results on low frequency estimates are relatively recent, compared to the longer history of the high frequency regime, and some of them are not yet optimal. The main result of this paper is a low frequency bound for the resolvent of the Laplace-Beltrami operator with sharp weight. The interest is twofold. On one hand, we obtain the same type of sharp inequality as on $ \Ra^n $ for a general class of manifolds which contains both $ \Ra^n $ with an asymptotically flat metric and the  class of scattering manifolds (see \cite{Melrose1,Melrose2}). On the other hand, in the spirit of the applications quoted above, our result can be used in the proof of global Strichartz estimates: it allows to handle in a fairly simple and intuitive fashion the phase space region which cannot be treated by semiclassical (or microlocal) techniques.

Let us describe more precisely our framework and our results.

In this paper we consider an asymptotically conical manifold $ ({\mathcal M},G) $, that is a Riemannian manifold isometric outside a compact subset to a product $ (R_0,+\infty) \times {\mathcal S} $, with $ ({\mathcal S},h_0 ) $ a closed Riemannian manifold, equipped with a metric approaching the conical metric $ dr^2 + r^2 h_0 $ as $ r \rightarrow \infty $. More precisely this means that for some compact, connected manifold with boundary $ {\mathcal K} \Subset {\mathcal M} $ and some $ R_0 >0 $, 
there is a diffeomorphism 
\begin{eqnarray}
\kappa :   {\mathcal M} \setminus {\mathcal K} \ni m \rightarrow  \big(r(m),\omega(m) \big) \in (R_0, + \infty) \times {\mathcal S} , \label{diffeokappa}
\end{eqnarray}
through which the metric reads 
\begin{eqnarray}
 G = \kappa^* \left( a(r)dr^2 + 2 r b (r)dr + r^2 h (r) \right)  , \label{formedeG}
\end{eqnarray} 
with $a \rightarrow 1$, $b \rightarrow 0$ and $h \rightarrow h_0$ as $ r \rightarrow \infty$ in the following sense: for each $ r > R_0 $, $ a (r) $ is a function on $ {\mathcal S} $, $ b (r) $ is a $1$-form on $ {\mathcal S}$ and $ h (r) $ is a Riemannian metric  on $ {\mathcal S} $, with $ a (\cdot) , b (\cdot) $ and $ h (\cdot) $ all depending smoothly on $r$  so  that,  for some $ \rho > 0 $,
\begin{eqnarray}
  || \partial_r^j (a (r) - 1) ||_{\Gamma^0 ({\mathcal S})} +  || \partial_r^j b (r)  ||_{\Gamma^1 ({\mathcal S})} +  || \partial_r^j (h (r) - h_0) ||_{\Gamma^2 ({\mathcal S})}  \lesssim r^{-j-\rho} ,   \label{hypothesedecroissance}
\end{eqnarray}
where, for $ k = 0,1,2 $, $ || \cdot ||_{\Gamma^k ({\mathcal S})} $ is any seminorm of the space of smooth sections of $ (T^* {\mathcal S})^{\otimes^k} $. In usual terms, this means that $ G $ is a long range perturbation of $ \kappa^* ( dr^2 + r^2 h_0 ) $ near infinity. 
%
%
In (\ref{diffeokappa}) $r$ is the first component of $ \kappa $. It defines a coordinate on $ {\mathcal M} \setminus {\mathcal K} $ taking its values in $ (R_0,\infty) $. We also assume that $ \kappa $ is  an homeomorphism between $ \overline{{\mathcal M} \setminus {\mathcal K}} $ and $ [R_0, \infty ) \times {\mathcal S} $. We may then assume without loss of generality that $r$ is a globally defined smooth function which is proper\footnote{{\it i.e.} $ r^{-1}( [r_1,r_2] ) $ is a compact subset of $ {\mathcal M} $ for all $ r_1 \leq r_2 $}, but which is a coordinate only near infinity. This  allows us to define the weights $ \scal{r}^{\mu} = (1+r^2)^{\mu/2}$ globally on $ {\mathcal M} $.


Our definition is more general than the one of scattering metrics \cite{Melrose1,Melrose2} and than the one used in \cite[Definition 1.4]{RoTa} where $ h $ has a polyhomogeneous expansion at infinity. It also covers the usual case of long range perturbations of the Euclidean metric as considered in \cite{BonyHafner,Boucletcanadian,BoucletCPDE}.

We will allow the possibility for $ {\mathcal M} $ to have a boundary. We thus introduce $ C_c^{\infty}({\mathcal M}) $, the set of smooth functions vanishing outside a compact set (these functions do not need to vanish on $ \partial {\mathcal M}$), and $ C_0^{\infty}({\mathcal M}) = C_c^{\infty} ({\mathcal M} \setminus \partial {\mathcal M}) $ the subset of those which also vanish near $ \partial M $. We let $ \hat{P} $ be the Friedrichs extension of $ - \Delta_G $ on $ C_0^{\infty}({\mathcal M}) $. It is self-adjoint on $ L^2 ({\mathcal M}) = L^2 ({\mathcal M}, d {\rm vol}_G) $. If $ {\mathcal M} $ has no boundary, it is the unique self-adjoint realization of $ - \Delta_G $ and if $ \partial {\mathcal M} $ is non empty it is the Dirichlet realization. 
We  assume that $ {\mathcal M} $ is connected to ensure that
\begin{eqnarray}
0 \ \mbox{is not an eigenvalue of} \ \hat{P} . \label{vpen0}
\end{eqnarray}
Our assumptions also imply that
\begin{eqnarray}
\chi (\hat{P}+i)^{-1} \ \mbox{is compact on} \ L^2 ({\mathcal M}) , \ \mbox{for all} \ \chi \in C_c^{\infty}({\mathcal M}) . \label{compacite}
\end{eqnarray}
We let $n = \mbox{dim} ( {\mathcal M} ) $ and assume everywhere that $ n \geq 3 $.

Our first main result is the following.


\begin{theo} \label{sharpweightresolvent} There exist $ \varepsilon_0 > 0 $ and $ C > 0 $ such that, for all $ z \in \Ca \setminus \Ra $ satisfying $ | \emph{Re}(z) | < \varepsilon_0 $,
$$ \big| \big| \scal{r}^{-1} ( \hat{P} - z )^{-1} \scal{r}^{-1} \big| \big|_{L^2 ({\mathcal M}) \rightarrow L^2 ({\mathcal M})} \leq C . $$
\end{theo}
 In \cite{GuHa1,Boucletcanadian,BonyHafner,GuHa2,BoucletCPDE,GuHaSi,RoTa}), uniform estimates on $  \scal{r}^{-s} ( \hat{P} - z )^{-1} \scal{r}^{-s} $ (for $ |\mbox{Re}(z)| $ small) were proved for  $s>1$.  The novelty of this result is that we use the weight $ \scal{r}^{-1} $ which is sharp.
  We  also cover more general manifolds than the ones considered in the aforementioned papers. To see the sharpness, we consider the flat Laplacian on $ \Ra^n $. If we could replace the weight $ \scal{r}^{-1} $ in Theorem \ref{sharpweightresolvent} by $ \scal{r}^{-s} $ for some $ s < 1 $, then by letting $z \rightarrow 0$,  we would obtain the $ L^2 $ boundedness of $ (1 + |x|)^{-s} \Delta^{-1} (1+|x|)^{-s} $ ($\Delta^{-1} $ being understood as the Fourier multiplier by $ -|\xi|^{-2} $) and then, by rescaling, we would have
  $$ \big| \big| (|x|+1)^{-s} \Delta^{-1} (|x|+1)^{-s} \big| \big|_{L^2 \rightarrow L^2} = \epsilon^{2-2s} \big| \big| (|x|+ \epsilon^{-1})^{-s} \Delta^{-1} (|x| + \epsilon^{-1})^{-s} \big| \big|_{L^2 \rightarrow L^2} \lesssim \epsilon^{2-2s} $$
  hence see that  $ (1+|x|)^{-s} \Delta^{-1} (1+|x|)^{-s} = 0 $ which is obviously wrong. 
  
This result is satisfactory for it answers the natural question of what the optimal weight is, but it also has  useful applications, in particular to the study of global in time Strichartz estimates which we describe below.

Our second main result is the following.

\begin{theo} \label{pseudorescaledestimates} Fix an integer $ N \geq 1 $ and a compact interval $ [E_1,E_2] \subset (0,\infty) $. There exist $ C > 0 $ and $ \epsilon_0 > 0 $ such that
\begin{eqnarray}
 \big| \big| \scal{\epsilon r}^{-N} \big( \epsilon^{-2} \hat{P} - Z \big)^{-N} \scal{\epsilon r}^{-N} \big| \big|_{L^2 ({\mathcal M}) \rightarrow L^2 ({\mathcal M})}
\leq C , \label{estimeesderniertheoreme}
\end{eqnarray}
for all $ Z \in \Ca \setminus \Ra $ such that $ \emph{Re}(Z) \in [E_1,E_2] $ and all $ \epsilon \in (0,\epsilon_0) $.  
\end{theo}

These estimates are  low frequency inequalities for they are equivalent to the spectrally localized versions
$$ \big| \big| \scal{\epsilon r}^{-N} \big( \epsilon^{-2} \hat{P} - Z \big)^{-N} \phi (\epsilon^{-2 } \hat{P} ) \scal{\epsilon r}^{-N} \big| \big|_{L^2 ({\mathcal M}) \rightarrow L^2 ({\mathcal M})}
\leq C $$
for any $ \phi \in C_0^{\infty} $ which is equal to $1$ near $ [ E_1 , E_2] $, that is when $ \hat{P} $ is spectrally localized near $ [\epsilon^2 E_1 , \epsilon^2 E_2] $. 
Let us remark that for the Laplacian on $ \Ra^n $,  Theorem \ref{pseudorescaledestimates} follows directly from the usual estimates on $ \scal{r}^{-N} (-\Delta - Z)^{-N} \scal{r}^{-N} $ by a simple rescaling argument.  Such a global rescaling argument is of course meaningless on a manifold, but Theorem \ref{pseudorescaledestimates} says that this scaling intuition remains correct. 

We will explain below to which extent Theorem \ref{pseudorescaledestimates} is complementary to Theorem \ref{sharpweightresolvent}. Before doing so, we record a last result which is a byproduct of our analysis  but which is also interesting on its own.
\begin{theo} \label{carreresolvante} Fix $ s \in (0,1/2) $. There exist $ \varepsilon_0 > 0 $ and $ C > 0 $ such that, for all $ z \in \Ca \setminus \Ra $ satisfying $ | \emph{Re}(z) |  < \varepsilon_0 $,
\begin{eqnarray}
 \big| \big| \scal{r}^{- 2 - s} ( \hat{P} - z )^{- 2} \scal{r}^{- 2 -s } \big| \big|_{L^2 ({\mathcal M}) \rightarrow L^2 ({\mathcal M})} \leq C |\emph{Re}(z)|^{s-1} . \label{nearlysharp}
\end{eqnarray} 
\end{theo}
The estimate (\ref{nearlysharp}) is nearly sharp with respect to $ |\mbox{Re}(z)|^{s-1} $ in dimension 3 (one can take exactly $ s = 1/2 $ in the asymptotically Euclidean case \cite{BoucletCPDE}, but this is not clear if $ {\mathcal S} \ne {\mathbb S}^2 $) but certainly not in higher dimensions (see \cite{BoucletCPDE} where we get better estimates  in higher dimensions in the asymptotically 
 Euclidean case). To get sharper estimates, one would need to use improved Hardy inequalities ({\it e.g.} improve Lemma \ref{pour2eme} to be able to consider higher order derivatives in higher dimensions). We did not consider this technical question since the main focus of this paper is on  Theorems \ref{sharpweightresolvent} and \ref{pseudorescaledestimates} for which Theorem \ref{carreresolvante} will  be essentially a tool.


We now discuss some  motivations and applications of Theorems \ref{sharpweightresolvent} and \ref{pseudorescaledestimates}. The first application is on the global smoothing effect.
\begin{coro} \label{corollairesmoothing} Assume that $ {\mathcal M} $ has no boundary and has no trapped geodesics. Then, there exists $  C>0 $ such that
$$ \int_{\Ra} || \scal{r}^{-1} (1+\hat{P})^{1/4} e^{it\hat{P}} u_0 ||_{L^2}^2 dt \leq C || u_0 ||_{L^2}^2. $$
\end{coro}
We state this result in the case of boundaryless manifolds only for simplicity. However, it extends to manifolds with boundary, under the non trapping condition for the generalized billard flow of Melrose-Sj\"ostrand (see \cite{MelroseSjostrand} and \cite{Burq} for related problems).

 On manifolds, the local in time version of this corollary is classical (see {\it e.g.} \cite{Doi}). We refer to \cite{Benartzi,Kato} for the global in time version in the flat case.
  Here we derive a global in time version with the sharp weight $ \scal{r}^{-1} $. According to the standard approach, Corollary \ref{corollairesmoothing} follows from the resolvent estimates
\begin{eqnarray}
 \big| \big| \scal{r}^{-1} (\hat{P}-\lambda \pm i \varepsilon)^{-1} \scal{r}^{-1} \big| \big|_{L^2 ({\mathcal M}) \rightarrow L^2 ({\mathcal M})} \leq C \scal{\lambda}^{-1/2}, \qquad \lambda \in \Ra , \ \varepsilon > 0 ,  \label{globalresolvent}
\end{eqnarray}
and the Kato theory of smooth operators  \cite{ReedSimon4}. The resolvent estimates (\ref{globalresolvent}) follow from \cite{Cardoso} at high energy, using the non trapping condition, from Theorem \ref{sharpweightresolvent} at low energy and, when $ \lambda $ belongs to any compact subset of $ (0,\infty) $, from the standard Mourre theory \cite{Mourre,JMP,FrHi} combined with the absence of embedded eigenvalues for $\hat{P}$ (see \cite{Ito}).

We next consider  Strichartz estimates, which is the main original motivation of this paper. The related results will appear in a forthcoming paper \cite{BoucletMizutani} but we explain below why Theorems \ref{sharpweightresolvent} and \ref{pseudorescaledestimates} are relevant to handle the contribution of low frequencies. In particular, we will see where  using the weight  $ \scal{r}^{-1} $ is crucial.

Assume that $ f \in C_0^{\infty} (0,\infty) $ is fixed and that we wish to prove spectrally localized (at low frequency) Strichartz estimates of the form
\begin{eqnarray}
 \left( \int_{\Ra} || f (\epsilon^{-2} \hat{P}) e^{it\hat{P}} u_0 ||_{L^{2^*}({\mathcal M})}^2 dt  \right)^{1/2} \leq C ||u_0 ||_{L^2({\mathcal M})} , \label{Strichartz}
\end{eqnarray} 
for some $ C $ independent of $ 0 < \epsilon < 1 $. Here $ 2^* = 2n / (n-2) $.
Choosing $ \chi \in C_0^{\infty} (\Ra) $ which is equal to $1$ near $ 0 $, we split $ u_{\epsilon}(t) := f (\epsilon^{-2} \hat{P}) e^{it\hat{P}} u_0 $ as the sum   $ \chi (\epsilon r) u_{\epsilon}(t) + (1 - \chi)(\epsilon r) u_{\epsilon}(t) $. By the homogeneous Sobolev  estimate, we have
\begin{eqnarray*}
   || \chi (\epsilon r) u_{\epsilon}(t) ||_{L^{2^*}({\mathcal M})} & \lesssim &  || \nabla_G \big( \chi (\epsilon r) u_{\epsilon}(t) \big)  ||_{L^2} \\
   & \lesssim & || \epsilon \scal{\epsilon r}^{-1} u_{\epsilon}(t) ||_{L^2} \\
   & \lesssim & || \scal{r}^{-1} u_{\epsilon} (t) ||_{L^2} .
\end{eqnarray*}   
The second estimate is not completely obvious (it is proved in \cite{BoucletMizutani} but is slightly too technical to be reproduced here). Heuristically, it follows from the fact that
\begin{eqnarray*}
 || \nabla_G \big( \chi (\epsilon r) u_{\epsilon}(t) \big)  ||_{L^2} & \lesssim & || \epsilon \chi^{\prime} (\epsilon r) u_{\epsilon}(t) \big)  ||_{L^2} + ||  \chi (\epsilon r) \nabla_G u_{\epsilon}(t) \big)  ||_{L^2} \\
 & \lesssim & || \epsilon \chi^{\prime} (\epsilon r) u_{\epsilon}(t) \big)  ||_{L^2} + || \epsilon  \chi (\epsilon r) (\epsilon^{-2} \hat{P} )^{1/2} u_{\epsilon}(t) \big)  ||_{L^2}
\end{eqnarray*}  
by formally replacing $ \nabla_G $ by $ \hat{P}^{1/2} (= \epsilon (\epsilon^{-2} \hat{P})^{1/2}) $. The point we want to emphasize here is that the homogeneity in $ \hat{P}^{1/2} $ combined with the spectral localization $ f (\epsilon^{-2} \hat{P}) $ allows to gain precisely one power of $ \epsilon $ which, combined with the localization $ \chi (\epsilon r) $ allows to gain a decay of order $ \scal{r}^{-1} $.  Then, by using Corollary \ref{corollairesmoothing} which implies that
\begin{eqnarray}
 \left( \int_{\Ra} || \scal{r}^{-1} f (\epsilon^{-2} \hat{P})  e^{it \hat{P}} u_0 ||_{L^2({\mathcal M})}^2 dt \right)^{1/2} \leq C || u_0 ||_{L^2({\mathcal M})} , 
 \nonumber \label{L2L2spacetime}
\end{eqnarray} 
 we get immediatly
$$ \left( \int_{\Ra} || \chi (\epsilon r) f (\epsilon^{-2} \hat{P}) e^{it\hat{P}} u_0 ||_{L^{2^*}({\mathcal M})}^2 dt  \right)^{1/2} \leq C || u_0 ||_{L^2} . $$
This argument uses crucially the sharp weight $ \scal{r}^{-1} $ of Theorem \ref{sharpweightresolvent}.
To get (\ref{Strichartz}), we next need  to consider $  (1- \chi) (\epsilon r) f (\epsilon^{-2} \hat{P}) e^{it\hat{P}} u_0 $. The corresponding analysis is far too long to be explained in detail in this paper, however we mention the following key idea. The interest of the localization $ (1- \chi) (\epsilon r) f (\epsilon^{-2} \hat{P}) $ is that, in the phase space, it corresponds to a region where $ |\xi|_G \sim \epsilon $ and $ r \gtrsim \epsilon^{-1} $ which shows,  on the basis of the uncertainty principle intuition, that we can use microlocal methods (typically the Isosaki-Kitada techniques as in \cite{Bouclet,Mizutani}). To implement this intuition concretely, we use that  the spatial localization $ 1-\chi (\epsilon r) $ allows to rescale our problem by a factor $ \epsilon^{-1} $ to work at frequency $1$, away from a compact subset of $ {\mathcal M} $. We can then  prove the relevant dispersion estimates by means of microlocal parametrices, whose remainders are controled thanks to the approriate resolvent (or propagation) estimates given by Theorem \ref{pseudorescaledestimates}. In summary, the interest of our results for Strichartz estimates are on one hand Theorem \ref{sharpweightresolvent}  which controls the uncertainty region $ \{ |\xi|_{G} \sim \epsilon , r \lesssim \epsilon^ {-1} \}$ and on the other hand Theorem  \ref{pseudorescaledestimates} which allows to study $ \{ |\xi|_{G} \sim \epsilon , r \gtrsim \epsilon^ {-1} \}$ via rescaling and microlocal techniques.

The paper is organized as follows. In Section \ref{section2}, we  prove Theorems \ref{sharpweightresolvent}, \ref{pseudorescaledestimates} and  \ref{carreresolvante} assuming  analogous inequalities on a model operator $ \hat{P}_T $. These estimates on $ \hat{P}_T $ follow from a suitable version of the Mourre Theory on a cone.
This theory  is described in  Section \ref{Mourreonacone} and applied in
 Section \ref{Mourreonaconesuite} to prove the results used in Section \ref{section2}.

\section{Proofs of the main results} \label{section2}

The basic idea to prove Theorems \ref{sharpweightresolvent}, \ref{pseudorescaledestimates} and  \ref{carreresolvante}  is to extrapolate $ {\mathcal M} \setminus {\mathcal K} \approx (R_0,\infty) \times {\mathcal S} $ into a pure cone $ {\mathcal M}_0= (0,\infty) \times {\mathcal S} $, where one can use a global scaling argument to reduce the proof to the case of estimates at frequency $1$. The contribution of $ {\mathcal K} $ is then treated by means of a compactness trick,  using that $ 0 $ is not an eigenvalue of $ \hat{P} $.

We record here the main steps of this analysis and then derive the proofs of Theorems \ref{sharpweightresolvent}, \ref{pseudorescaledestimates} and  \ref{carreresolvante}. 

\begin{prop} \label{reductiongeometriqueprop} We may assume that the diffeomorphism $ \kappa $ in (\ref{diffeokappa}) is such that 
\begin{eqnarray}
 d {\rm vol}_G = \kappa^* \big( r^{n-1} dr d {\rm vol}_{h_0} \big) , \qquad  \mbox{on} \  {\mathcal M} \setminus {\mathcal K}  . \label{formedesniterevue}
\end{eqnarray} 
 More precisely, we can find a new $ \tilde{\mathcal K} \Subset {\mathcal M} $, a new diffeomorphism $ \tilde{\kappa} $ and a new proper function $ \tilde{r} : {\mathcal M} \rightarrow [ 0 , \infty ) $, such that  (\ref{diffeokappa}), (\ref{formedeG}), (\ref{hypothesedecroissance})  hold with $(\tilde{\kappa} , \tilde{r} ) $,  such that  $ \scal{r} / \scal{\tilde{r}} $ is bounded from above and below on $ {\mathcal M} $ by positive constants, and such that $ d {\rm vol}_G = \tilde{\kappa}^* \big( \tilde{r}^{n-1} d \tilde{r} d {\rm vol}_{h_0} \big)  $ on $ {\mathcal M} \setminus \tilde{\mathcal K} $.
\end{prop}

\noindent {\it Proof.} See Section \ref{Mourreonaconesuite}.

\bigskip

The interest of working with the density $ r^{n-1}dr d {\rm vol}_{h_0} $  is that, on the cone $ (0,\infty) \times {\mathcal S} $, the group $ (e^{itA})_{t \in \Ra} $ of $ L^2 $ scalings (see (\ref{dilations}))
is unitary on $ L^2 \big( (0,\infty) \times {\mathcal S} , r^{n-1} dr d {\rm vol}_{h_0} \big)  $. This guarantees that both $A$ and a suitable extrapolation $ \hat{P}_T $ of $ \hat{P} $ on $ (0,\infty) \times {\mathcal S} $ (see Proposition \ref{parametrixe} below) are  self-adjoint with respect to the same measure.

We next record useful results related to the Hardy inequality. We define the operator
$$  \hat{P}^{-1/2} := f (\hat{P}) , \qquad f (\lambda) = {\mathds 1}_{(0,+\infty)}(\lambda) \lambda^{-1/2} , $$
 by means of the spectral theorem (see {\it e.g.} \cite[p. 263]{ReedSimon1}). It is an unbounded self-adjoint operator and it is a routine  to check that 
 it maps its domain into the domain of $ \hat{P}^{1/2} $. Moreover, we have 
\begin{eqnarray}
 \hat{P}^{1/2} \hat{P}^{-1/2} = I , \qquad  \mbox{on Dom}(\hat{P}^{-1/2}) , \label{algebreracine}
\end{eqnarray}
which is a  consequence of the spectral theorem and the property (\ref{vpen0}).
\begin{prop} \label{Hardyglobal} \begin{enumerate} \item{ There is a constant $ C $ such that, for all $ u \in \emph{Dom}\big(\hat{P}^{1/2} \big) $,
\begin{eqnarray}
 || \scal{r}^{-1} u ||_{L^2 ({\mathcal M})} \leq C \big| \big| \hat{P}^{1/2} u \big| \big|_{L^2 ({\mathcal M})} . \label{Hardylabel}
\end{eqnarray} }
\item{The operator $ \scal{r}^{-1} $ maps $ L^2({\mathcal M}) $ into $ \emph{Dom}(\hat{P}^{-1/2}) $ and $ \hat{P}^{-1/2} \scal{r}^{-1} $ is bounded on $ L^2 ({\mathcal M}) $.}
\item{For all real number $ s > 1 $, $ \hat{P}^{-1/2} \scal{r}^{-s} $ is compact on $ L^2 ({\mathcal M}) $.}
\end{enumerate}
\end{prop}

\noindent {\it Proof.} 1. Since $ \hat{P} $ is the Friedrichs extension of $ - \Delta_G $, $ C_0^{\infty} ({\mathcal M}) $  is dense in $ \mbox{Dom}\big( \hat{P}^{1/2} \big) $ for the graph norm so we may assume that $ u \in C_0^{\infty}({\mathcal M}) $. Fix $ \chi = \chi (r) $ a smooth function which is equal to $ 1 $ near $ {\mathcal K} $ and vanishes for $ r \gg 1 $. Then
$$ || \scal{r}^{-1} u ||_{L^2 ({\mathcal M})} \leq  || \chi u ||_{L^2 ({\mathcal M})} + || \scal{r}^{-1} (1-\chi) u ||_{L^2 ({\mathcal M})} . $$
From the Poincar\'e inequality on compact manifold with boundary containing $ \mbox{supp} (\chi) $, we get
$$ || \chi u ||_{L^2 ({\mathcal M})} \lesssim || \nabla_G (\chi u) ||_{L^2({\mathcal M})} \lesssim || \nabla_G u  ||_{L^2 ({\mathcal M})} + || \chi^{\prime} (r) u ||_{L^2 ({\mathcal M})} . $$
On the other hand, using the Hardy inequality on $ (0,\infty) \times {\mathcal S} $ (see  (\ref{Hardy1}) in Section \ref{Mourreonacone}) we also have
$$ ||\scal{r}^{-1} (1- \chi) u ||_{L^2 ({\mathcal M})} \lesssim || \partial_r ((1-\chi) u) ||_{L^2({\mathcal M})} \lesssim || \nabla_G u  ||_{L^2 ({\mathcal M})} + || \chi^{\prime} (r) u ||_{L^2 ({\mathcal M})} . $$
To complete the proof it suffices to observe that, if $ \chi^{\prime}(r) $ is supported in $ \{r_1 < r < r_2 \} $ with $r_1 > R_0$,
\begin{eqnarray}
 || \chi^{\prime} (r) u ||_{L^2 ({\mathcal M})} \lesssim ||  u ||_{L^2 (r_1 < r < r_2)} \lesssim || \nabla_G u ||_{L^2 ({\mathcal M})} . \label{transitionab}
\end{eqnarray} 
Indeed, for all $ r > R_0 $ and all $ \omega \in {\mathcal S} $, we have
$$ | \kappa^* u (r,\omega) | = \left| \int_r^{\infty} \partial_s (\kappa^ *u) (s, \omega) ds \right| \leq \frac{r^{\frac{2-n}{2}}}{(n-2)^{1/2}}\left( \int_r^{\infty} | \partial_s (\kappa^ *u) (s,\omega)|^2 s^{n-1} ds \right)^{1/2} , $$
using the Cauchy-Schwartz inequality. Squaring and integrating over $ [r_1,r_2] \times{\mathcal S} $ with respect to $ r^{n-1} dr d {\rm vol}_{h_0} (\omega) $, we get (\ref{transitionab}).
Since $ || \nabla_G u ||_{L^2 ({\mathcal M})} = \big| \big| \hat{P}^{1/2} u \big| \big|_{L^2 ({\mathcal M})} $, (\ref{Hardylabel}) follows. 

\noindent 2. For all $ u \in \mbox{Dom} (\hat{P}^{-1/2}) $ and $ v \in L^2 ({\mathcal M}) $, we have
$$ \big| \big( \scal{r}^{-1} v , \hat{P}^{-1/2} u \big) \big| = \big| \big(  v , \scal{r}^{-1} \hat{P}^{-1/2} u \big) \big| \leq C || v ||_{L^2 ({\mathcal M})} || u ||_{L^2 ({\mathcal M})} , $$
using (\ref{algebreracine}) and (\ref{Hardylabel}). This implies both that $ \scal{r}^{-1}v $ belongs to $ \mbox{Dom} \big( (\hat{P}^{-1/2})^* \big) = \mbox{Dom}(\hat{P}^{-1/2}) $ and that $ \hat{P}^{-1/2} \scal{r}^{-1}  $ is bounded. 

\noindent 3. Fix $ \Phi \in C_0^{\infty}(\Ra) $ which is equal to $1$ near zero. Then, using the spectral theorem, we have the decomposition
$$ \hat{P}^{-1/2} \scal{r}^{-1} = \Phi(\hat{P}) \hat{P}^{-1/2} \scal{r}^{-s} + (1 - \Phi)(\hat{P}) \hat{P}^{-1/2} \scal{r}^{-s} . $$
The second term in the right-hand side is compact since (\ref{compacite}) implies that $ f (\hat{P})g(r) $ is compact whenever $f$ and $g$ are continuous on $ \Ra $ and vanish at infinity. We then rewrite the first term of the right-hand side as $ \hat{P}^{-1/2} \scal{r}^{-1} \big( \scal{r} \Phi (\hat{P}) \scal{r}^{-s} \big) $ so the result will follow from the compactness of $ \scal{r} \Phi (\hat{P}) \scal{r}^{-s} $ proved as follows. By the Helffer-Sj\"ostrand formula \cite{DiSj}, we have
\begin{eqnarray}
 \scal{r} \Phi (\hat{P}) \scal{r}^{-s} = \int_{\Ca} \overline{\partial} \widetilde{\Phi}(z) \scal{r} (\hat{P}-z)^{-1} \scal{r}^{-s} L (dz) , \label{HelfferSjostrand}
\end{eqnarray}
where $ \widetilde{\Phi} \in C_c^{\infty}(\Ca) $ is an almost analytic extension of $ \Phi $, {\it i.e.} $ \overline{\partial} \widetilde{\Phi}(z) = {\mathcal O} (\mbox{Im}(z)^{\infty}) $, and
$$ \scal{r} (\hat{P}-z)^{-1} \scal{r}^{-s} = (\hat{P}-z)^{-1} \scal{r}^{-(s-1)} + (\hat{P}-z)^{-1} [\hat{P},\scal{r}] (\hat{P}-z)^{-1} \scal{r}^{-s} . $$
Using on one hand that $ (\hat{P}-z)^{-1} [\hat{P},\scal{r}] $ has a bounded closure on $ L^2 ({\mathcal M}) $ with norm $ {\mathcal O}(\mbox{Im}(z)^{-1}) $ on  $ \mbox{supp}(\widetilde{\Phi}) $, and  on the other hand that, for $ \mu = s $ or $ s-1 $,  
$$ (\hat{P}-z)^{-1} \scal{r}^{-\mu} = (\hat{P}+i)(\hat{P}-z)^{-1} (\hat{P}+i)^{-1} \scal{r}^{-s} , $$
we obtain easily the compactness of (\ref{HelfferSjostrand}) from (\ref{compacite}). This completes the proof.
\finpreuve

\bigskip

\noindent {\bf Remark.} In the sequel we shall use freely the fact that, for any real number $ M $, operators of the form $ \scal{r}^{M} (\hat{P}+1)^{-1} \scal{r}^{-M} $ or $ \scal{r}^M \Phi (\hat{P}) \scal{r}^{-M} $, with $ \Phi \in C_0^{\infty}(\Ra) $, are bounded on $ L^2 ({\mathcal M}) $ and even map $ L^2 ({\mathcal M}) $ into $ \mbox{Dom}(\hat{P}^{1/2}) $. This is basically well known and follows from the same type of standard argument as the ones used in the proof of the item 3 above.

\bigskip

In the next proposition, we summarize the main technical results of this paper, which deal with resolvent estimates of model operators on
\begin{eqnarray}
 {\mathcal M}_0 := (0,+\infty) \times {\mathcal S}  , \label{Mzero}
\end{eqnarray}
equipped with the conical volume density $ r^{n-1} dr d {\rm vol}_{h_0} $. We will set everywhere 
$$ L^2 ({\mathcal M}_0) := L^2 \big( (0,+\infty) \times {\mathcal S} , r^{n-1} dr d {\rm vol}_{h_0} \big) . $$
In the sequel, when nothing is specified, $ || \cdot || $ will denote both norms $ || \cdot ||_{L^2 ({\mathcal M}_0) \rightarrow L^2 ({\mathcal M}_0)} $ and $ || \cdot ||_{L^2 ({\mathcal M})
\rightarrow L^2 ({\mathcal M})} $. This won't cause any ambiguity in practice but will simplify the notation.

We will also use the standard notation  $ \kappa^* $ for the  composition with $ \kappa $ and $ \kappa_* $ for the composition with $ \kappa^{-1} $.
\begin{prop} \label{parametrixe} Fix $ N \geq 1 $. There exists a self-adjoint operator $ \hat{P}_T $  on $ L^2 ({\mathcal M}_0) $ such that, 
for some $ R \gg 1 $, we have
$$ \kappa^* \hat{P}_T = \hat{P} \kappa^* , \qquad \mbox{on} \ (R,\infty) \times {\mathcal S} , $$
and such that $ \hat{P}_T $ satisfies the following resolvent estimates:
\begin{enumerate}
\item{there exists $ C  $ such that 
\begin{eqnarray}
\big| \big| r^{-1} \big(  \hat{P}_T - z \big )^{-1}  r^{-1} \big| \big|_{ L^2 ({\mathcal M}_0) \rightarrow L^2 ({\mathcal M}_0) } \leq C , 
\label{1ertheoreme}
\end{eqnarray}
for all $z\in \Ca \setminus \Ra$ such that $ |\emph{Re}(z)| \leq 1 $.}
\item{ For all $ [E_1,E_2] \Subset (0, + \infty) $ and all $ 0 \leq k \leq N-1 $, there exists $ C_k  $ such that
\begin{eqnarray}
\big| \big| (\epsilon r)^{-1} \scal{\epsilon r}^{-k} \big( \epsilon^{-2} \hat{P}_T - Z \big )^{-1-k} \scal{\epsilon r}^{-k} (\epsilon r)^{-1} \big| \big|_{ L^2 ({\mathcal M}_0) \rightarrow L^2 ({\mathcal M}_0) } \leq C_k , \label{rescalezero}
\end{eqnarray}
for all $ Z \in \Ca \setminus \Ra $ such that $ \emph{Re}(Z) \in [E_1,E_2] $ and all $ \epsilon \in (0,1] $.}
\item{ For all $ s \in (0,1/2) $, there exists $ C_s  $ such that
\begin{eqnarray}
\big| \big| \scal{ r}^{-2-s} (\hat{P}_T - z)^{-2} \scal{ r}^{-2 -s}  \big| \big|_{L^2 ({\mathcal M}_0) \rightarrow L^2 ({\mathcal M}_0)} \leq C_s |\emph{Re}(z)|^{s-1},  \label{pourcompacite}
\end{eqnarray}
for all $ z \in \Ca \setminus \Ra $ such that $ 0 < | \emph{Re}(z) | \leq 1 $.}
\end{enumerate} 
\end{prop}

\noindent {\it Proof.} See Section \ref{Mourreonaconesuite}. We only mention here that the operator $ \hat{P}_T $ will be constructed as the self-adjoint realization of an operator in divergence form given by (\ref{pourestimationpedagogie}).

\bigskip

Note that the first item of Proposition \ref{parametrixe} means essentially that $ \hat{P} $ and $ \hat{P}_T $ coincide close to infinity, though they are not defined on the same manifold. 

We then want to compare the resolvent of $ \hat{P} $ with the resolvent of $ \hat{P}_T $. Since they are not defined on the same space, we will use the following identification operators,
\begin{eqnarray*}
 J_0 & = & \chi_{{\mathcal M} \setminus {\mathcal K}} \circ \kappa^* \circ  \varrho_{(R_0,\infty) \times {\mathcal S}} : L^2 ({\mathcal M}_0 ) \rightarrow L^2 ({\mathcal M}) , \\
 J & = & \chi_{(R_0,\infty) \times {\mathcal S}} \circ \kappa_* \circ \varrho_{{\mathcal M} \setminus {\mathcal K}} : L^2 ({\mathcal M}) \rightarrow L^2 ({\mathcal M}_0)
\end{eqnarray*}
 Everywhere $ \varrho_{\Omega} $ stands for the restriction operator to $ \Omega $. It is straightforward to check that  $ J_0 J : L^2 ({\mathcal M}) \rightarrow L^2({\mathcal M}) $ is the multiplication operator
\begin{eqnarray}
 J_0 J = {\mathds 1}_{{\mathcal M} \setminus {\mathcal K}} . \label{J0J}
\end{eqnarray}
\begin{prop}[Generalized resolvent identity] \label{propositionasimplifer} Let $ \psi \in C^{\infty} ( {\mathcal M} ) $ depend only on  $r$, be supported  on $ \{ r  > R \} $ with $ R $ as in Proposition \ref{parametrixe}, and $ \psi \equiv 1 $ near infinity. Then, for all $ z \in \Ca \setminus \Ra $,
\begin{eqnarray}
(\hat{P}-z)^{-1} \psi J_0 & = & \psi J_0 (\hat{P}_T-z)^{-1} - (\hat{P}-z)^{-1} [\hat{P},\psi] J_0 (\hat{P}_T-z)^{-1} , \label{reso1} \\
J \psi (\hat{P}-z)^{-1} & = & (\hat{P}_T-z)^{-1} J \psi + (\hat{P}_T-z)^{-1} J [\hat{P},\psi] (\hat{P}-z)^{-1} . \label{reso2}
\end{eqnarray}
As a consequence, we also have
\begin{eqnarray}
(\hat{P}-z)^{-1} \psi  & = & \psi J_0 (\hat{P}_T-z)^{-1} J - (\hat{P}-z)^{-1} [\hat{P},\psi] J_0 (\hat{P}_T-z)^{-1} J , \label{reso11} \\
 \psi (\hat{P}-z)^{-1} & = & J_0 (\hat{P}_T-z)^{-1} J \psi + J_0 (\hat{P}_T-z)^{-1} J [\hat{P},\psi] (\hat{P}-z)^{-1} . \label{reso22}
\end{eqnarray}
\end{prop}

\noindent {\it Proof.}  By the first item of Proposition \ref{parametrixe} and the support property of $ \psi $, we have $ \psi \hat{P} J_0 = \psi J_0 \hat{P}_T $ and therefore
$$ (\hat{P} - z ) \psi J_0 (\hat{P}_T-z)^{-1} = [\hat{P},\psi] J_0 (\hat{P}_T-z)^{-1} + \psi J_0 . $$
After composition to the left with $ (\hat{P}-z)^{-1} $ we obtain (\ref{reso1}). The proof of (\ref{reso2}) is similar. The identities (\ref{reso11}) and (\ref{reso22}) follow from (\ref{reso1}) and (\ref{reso2}) respectively, combined with (\ref{J0J}) and the fact that $ \psi {\mathds 1}_{{\mathcal M} \setminus {\mathcal K}} = \psi $. \finpreuve

\bigskip



In the sequel, we fix a real number $ M > \max (N,3) $, $ N$ being the integer fixed in Proposition \ref{parametrixe}, and introduce the weighted resolvent
$$ R_M (z) = \scal{r}^{-M} (\hat{P} - z)^{-1} \scal{r}^{-M} . $$
We will first obtain estimates on $ R_M (z) $ and its derivatives. With those estimates at hand, it will be fairly easy to derive estimates with sharper weights.
To prove estimates on $ R_M (z) $ we will use the following proposition.

\begin{prop} \label{propvide} One can find $ \varepsilon_0 > 0 $ and bounded operators
$ S_1,S_2,T_1,T_2 $ and $ F_M (\cdot) , S(\cdot), T (\cdot) $ depending holomorphically on $ z \in \{ |\emph{Re}(z)| < \varepsilon_0 \} \cap \{ \emph{Im}(z) \ne 0 \} $ such that for all $ z  $ in this set, we have
\begin{eqnarray}
 R_M (z) = F_M (z) + S_1 R_M (z) T_1 + T_2 R_M (z) S_2 + T (z) R_M (z) S (z) ,  \label{neumann}
\end{eqnarray}
and
\begin{eqnarray}
 \big| \big| S_1 \big| \big| \ \big| \big| T_1 \big| \big| + \big| \big| S_2 \big| \big| \ \big| \big| T_2 \big| \big| + \big| \big| S (z) \big| \big| \ \big| \big| T (z) \big| \big| \leq 3 /4 , \label{smallness} 
\end{eqnarray}
and such that, for all $ s \in (0,1/2) $ and all $ 1 \leq k \leq N $,
\begin{eqnarray}
|| F_M (z) || & \leq & C, \qquad \qquad \qquad \qquad \ \ |\emph{Re}(z)| < \varepsilon_0 , \label{borne0} \\
|| \partial_z F_M (z) || + || \partial_z S (z) || + || \partial_z T (z) || & \leq & C_s |\emph{Re}(z)|^{s-1} , \qquad  0 \ne |\emph{Re}(z)| < \varepsilon_0 , \label{borne1} \\
|| \partial_z^k F_M (z) || + || \partial_z^k S (z) || + || \partial_z^k T (z) || & \leq & C_k \emph{Re}(z)^{-k} , \qquad \ \ \ \  0 < \emph{Re}(z) < \varepsilon_0 . \label{borne2}
\end{eqnarray}
\end{prop}

This proposition is based on Proposition \ref{propositionasimplifer}. Before proving it, we need to establish several intermediate lemmas. Consider three functions $ \psi_1 , \psi_2 $ and $ \varphi $ in $ C^{\infty}({\mathcal M}) $ supported in $ \{ r > R \} $, depending only on $r$, equal to $ 1 $ near infinity and such that 
$$ \psi_2 \equiv 1 \ \ \mbox{near supp}( \psi_1 ), \qquad \psi_1 \equiv 1 \ \ \mbox{near supp}(\varphi) . $$  
By  Proposition \ref{propositionasimplifer}
and the easily verified fact that $ J [\hat{P},\psi_1] J_0 = [\hat{P}_T,\psi_1] $ (we identify $ \psi_1  $ and $ \kappa_* \psi_1 $ in the obvious fashion since $ \psi_1 $ depends only on $r$), we obtain
\begin{eqnarray}
\varphi (\hat{P}-z)^{-1} \varphi & = & \varphi J_0 (\hat{P}_T -z)^{-1} J \varphi - \varphi J_0 (\hat{P}_T - z)^{-1} \big[ \hat{P}_T,\psi_1 \big] (\hat{P}_T - z)^{-1} J \varphi \nonumber \\
& & - \ \varphi J_0 (\hat{P}_T-z)^{-1} J [\hat{P},\psi_2] (\hat{P}-z)^{-1} [\hat{P},\psi_1] J_0 (\hat{P}_T-z)^{-1} J \varphi . \label{secondline}
\end{eqnarray}
The interest of this formula is that  $ [\hat{P},\psi_1] $ and $ [\hat{P},\psi_2] $ (as well as $ [\hat{P}_T,\psi_1] $) have compactly supported coefficients.
The smallness condition (\ref{smallness}) will be a consequence of the following lemma.

\begin{lemm} \label{lemmeacitermieux} Let $ \Phi \in C_0^{\infty}(\Ra) $ be equal to $1$ near $0$.
For all $ \delta > 0 $, we can choose $  \varepsilon > 0  $ and $ \nu >0 $ such that
\begin{eqnarray}
\big| \big| \scal{r}^M (1-\varphi)  \Phi (\hat{P}/\varepsilon) \big| \big| \leq \delta , \label{compactfacile}
\end{eqnarray}
and, for all $ z \in \Ca \setminus \Ra $ satisfying $ |\emph{Re}(z)|< \nu $,
\begin{eqnarray}
\big| \big|  \scal{r}^{-M} J_0 (\hat{P}_T - z)^{-1} J \varphi \Phi (\hat{P}/ \varepsilon) \scal{r}^{-M} \big| \big| \leq \delta  .
\label{compactHardy}
\end{eqnarray}
\end{lemm}

\noindent {\it Proof.} The first inequality is standard. We recall the proof for completeness. If $ \varepsilon $ is small enough, we have $ \Phi(\lambda/\varepsilon) = \Phi (\lambda) \Phi (\lambda/\varepsilon) $ so that
$$ \scal{r}^M (1-\varphi)  \Phi (\hat{P}/\varepsilon) = \left( \scal{r}^M (1-\varphi)  \Phi (\hat{P}) \right) \Phi (\hat{P}/\varepsilon) . $$
By (\ref{compacite}) and the compact support of $ 1-\varphi $, the first term in the right-hand side is a compact operator. On the other hand, the property (\ref{vpen0}) implies that $ \Phi (\hat{P}/\varepsilon) \rightarrow 0  $ in the strong sense as $ \varepsilon \rightarrow 0 $, hence in operator norm when composed with a compact operator. This yields  (\ref{compactfacile}). Let us prove
(\ref{compactHardy}). We set $ y = {\rm Im} (z) $ and split the resolvent as $$  (\hat{P}_T - z)^{-1} = (\hat{P}_T - i y)^{-1} + \left( (\hat{P}_T - z)^{-1} - (\hat{P}_T - i y)^{-1} \right) . $$
 We start with the contribution of the  first term which, using Proposition \ref{Hardyglobal}, we write  as
$$ \scal{r}^{-M} (\hat{P}_T - i y)^{-1} J \varphi \Phi (\hat{P}/\varepsilon) \scal{r}^{-M} = \big( \scal{r}^{-M} (\hat{P}_T - i y)^{-1} J \varphi  \hat{P}^{1/2} \big) \big( \Phi (\hat{P}/\varepsilon) \hat{P}^{-1/2} \scal{r}^{-M} \big) . $$
By Proposition \ref{Hardyglobal},  $ \hat{P}^{-1/2} \scal{r}^{-M} $ is compact, hence the second bracket in the right-hand side goes to zero in operator norm as $ \varepsilon \rightarrow 0 $. Therefore, it suffices to obtain a uniform estimate on the first term. To do so, we use the form of $ \hat{P}_T $ which is (the Friedrichs extension of) the elliptic operator in divergence form (\ref{pourestimationpedagogie}). Using the form of $ \hat{P}_T $, the Hardy inequality (\ref{Hardy1}) on $ {\mathcal M}_0 $ and the spectral theorem, we have
\begin{eqnarray*}
\big| \big| \scal{r}^{-M} (\hat{P}_T - i y)^{-1} J \varphi  \hat{P}^{1/2} \big| \big| & = & \big| \big| \hat{P}^{1/2} \varphi J^* (\hat{P}_T + i y)^{-1} \scal{r}^{-M}  \big| \big| \\
& \lesssim &  \big| \big| \nabla_G \varphi J^* (\hat{P}_T + i y)^{-1} \scal{r}^{-1}  \big| \big| \\
& \lesssim &   \big| \big| \hat{P}_T^{1/2} (\hat{P}_T + i y)^{-1} \scal{r}^{-1}  \big| \big| +  \big| \big| \scal{r}^{-1} (\hat{P}_T + i y)^{-1} \scal{r}^{-1}  \big| \big| \\
& \lesssim &  \big| \big| \hat{P}_T^{1/2} (\hat{P}_T + i y)^{-1} \hat{P}_T^{1/2}  \big| \big| \lesssim 1 .
\end{eqnarray*}
 We can therefore fix $ \varepsilon $ such that
$$ \big| \big| \scal{r}^{-M} (\hat{P}_T - i y)^{-1} J \varphi \Phi (\hat{P}/\varepsilon) \scal{r}^{-M} \big| \big| \leq \delta / 2 . $$
Then, by writing
$$ (\hat{P}_T - z)^{-1} - (\hat{P}_T - i y)^{-1} = \int_0^{{\rm Re}(z)} (\hat{P}_T -x-iy)^{-2}dx $$
and using (\ref{pourcompacite}), say with $ s = 1/4 $, we obtain
\begin{eqnarray*}
 \big| \big| \scal{r}^{-M} \left( (\hat{P}_T - z)^{-1} - (\hat{P}_T - i y)^{-1} \right) J \scal{r}^{-M} \big( \scal{r}^M \varphi \Phi (\hat{P}/\varepsilon) \scal{r}^{-M} \big) \big| \big| \leq  C_{\varepsilon} |\mbox{Re}(z)|^{1/4} .
\end{eqnarray*}
Here we also used that $ J $ preserves the decay in $ r $. The above norm can therefore be made smaller than $ \delta/2 $ if $ \mbox{Re}(z) $ is small enough and the result follows. \finpreuve

\bigskip

\begin{lemm} \label{souslemme1} Fix two integers $ k , M \geq 0 $. Then there exists $ C > 0$ such that
\begin{eqnarray*}
\big| \big| \partial_z^k \big( w (\hat{P}_T-z)^{-1} [\hat{P}_T, \psi_1] (\hat{P}_T-z)^{-1} w \big) \big| \big| \ \leq  
 \      C  \sum_{k_1 = 0}^k || w (\hat{P}_T-z)^{-1-k_1} \scal{r}^{-M} || \ +  \qquad \\
\qquad \ C \sum_{k_1 + k_2 = k} \big| \big| w (\hat{P}_T-z)^{-1-k_1} \scal{r}^{-M} \big| \big| \  \big| \big| \scal{r}^{-M} (\hat{P}_T-z)^{-1-k_2} w \big| \big| ,
\end{eqnarray*} 
for all $ z \in \Ca \setminus \Ra $ such that $ |\emph{Re}(z)| \leq 1 $ and all  weight $w$ such that $ ||w||_{L^{\infty}} \leq 1 $.
\end{lemm}

\noindent {\it Proof.} Fix $ \Phi \in C_0^{\infty} (\Ra) $ such that $ \Phi \equiv 1 $ near $ [-1,1] $ and write
$$ [\hat{P}_T,\psi_1] (\hat{P}_T-z)^{-1} = [\hat{P}_T,\psi_1] \Phi (\hat{P}_T) (\hat{P}_T-z)^{-1} + [\hat{P}_T,\psi_1] (1-\Phi)(\hat{P}_T) (\hat{P}_T-z)^{-1} . $$
The last term can be written $  [\hat{P}_T,\psi_1] (\hat{P}_T+i)^{-1} (1-\Phi)(\hat{P}_T) (\hat{P}_T+i) (\hat{P}_T-z)^{-1} $, which is holomorphic and bounded (uniformly in $z$) in the strip $ \{| \mbox{Re}(z) | \leq 1 \} $. Note also that, by the compact support of $ [\hat{P}_T,\psi_1] $ the same property holds for $ \scal{r}^{M} [\hat{P}_T,\psi_1]  (1-\Phi)(\hat{P}_T) (\hat{P}_T-z)^{-1} $. On the other hand, $ [ \hat{P}_T,\psi_1] \Phi (\hat{P}_T) $ is bounded and, by a Helffer-Sj\"ostrand formula argument (see the remark after Proposition \ref{Hardyglobal}), $ \scal{r}^M  [ \hat{P}_T,\psi_1] \Phi (\hat{P}_T) \scal{r}^M $ is bounded too. All this implies that
$$ (\hat{P}_T-z)^{-1}[\hat{P}_T,\psi_1] (\hat{P}_T-z)^{-1} = (\hat{P}_T-z)^{-1} \scal{r}^{-M} B_M \scal{r}^{-M} (\hat{P}_T-z)^{-1} + (\hat{P}_T-z)^{-1} \scal{r}^{-M} O_M (z) , $$
for some bounded operators $B_M $ and $ O_M (z) $, $ O_M (\cdot)$ being holomorphic and uniformly bounded in the strip $ \{ |\mbox{Re}(z)| \leq 1 \} $. The result then follows easily. \finpreuve

\bigskip

\begin{lemm} \label{souslemme2} Fix $ 0 \leq k \leq N -1 $. Then 
\begin{eqnarray}
\big| \big| \scal{r}^{-M} ( \hat{P}_T - z)^{-1-k} \scal{r}^{-M}  \big| \big|
& \lesssim & \emph{Re}(z)^{-k} \label{showing} , \\
\big| \big| \scal{r}^{-M} \partial_z^k \big(  ( \hat{P}_T - z)^{-1}[\hat{P}_T,\psi_1] (\hat{P}_T-z)^{-1} \big) \scal{r}^{-M} \big| \big|
& \lesssim & \emph{Re}(z)^{-k} , \label{showingbis}
\end{eqnarray}
for all $ z \in \Ca \setminus \Ra $ such that $ \emph{Re}(z)> 0 $ is small enough
\end{lemm}

\noindent {\it Proof.} We rewrite $ z = \epsilon^2 Z $ with $ \mbox{Re}(Z) $ in a compact subset of $ (0,\infty) $ and $ \epsilon \in (0,1] $. 
Then, by using (\ref{rescalezero}) and the fact that
$$ (\epsilon r)^{-1} \scal{\epsilon r}^{-k} ( \epsilon^{-2} \hat{P}_T - Z)^{-1-k} \scal{\epsilon r}^{-k} (\epsilon r)^{-1}  = \epsilon^{2k} r^{-1} \scal{\epsilon r}^{-k}  (\hat{P}_T - \epsilon^2 Z)^{-1-k}   \scal{\epsilon r}^{-k} r^{-1} , $$
whose left-hand side is bounded by (\ref{rescalezero}), we get (\ref{showing}) once observed that
$$ \scal{r}^{-M} = \big( \scal{r}^{-M} r \scal{\epsilon r}^k \big) r^{-1} \scal{\epsilon r}^{-k} , $$
where the bracket is bounded uniformly in $ \epsilon $. The estimate (\ref{showingbis}) follows from Lemma \ref{souslemme1} and (\ref{showing}). \finpreuve

\bigskip

\noindent {\bf Proof of Proposition \ref{propvide}.} 
 We use the spatial cutoffs $ \varphi , \psi_1 , \psi_2 $ introduced after Proposition \ref{propvide}. We write
\begin{eqnarray}
 (\hat{P}-z)^{-1} = \varphi (\hat{P}-z)^{-1} \varphi + \varphi (\hat{P}-z)^{-1}(1-\varphi) + (1-\varphi) (\hat{P}-z)^{-1} , \label{partitiondelunitespatiale}
\end{eqnarray}
where we recall that $ 1 - \varphi $ is compactly supported on $ {\mathcal M} $. We next consider three functions $ \Phi , \phi_1 , \phi_2 \in C_0^{\infty} (\Ra) $ equal to $1$ near $ 0 $, which we will use as spectral cutoffs. We fix $ \Phi $ arbitrarily but we will choose $ \phi_1 , \phi_2 $ below. We use $ \Phi $ to rewrite the second line of (\ref{secondline}) as 
the sum of
\begin{eqnarray}
-  \varphi J_0 (\hat{P}_T-z)^{-1} J \left\{ [\hat{P},\psi_2] \Phi (\hat{P}) (\hat{P}-z)^{-1} \Phi (\hat{P}) [\hat{P},\psi_1] \right\} J_0 (\hat{P}_T-z)^{-1} J \varphi
\end{eqnarray}
and
\begin{eqnarray}
-  \varphi J_0 (\hat{P}_T-z)^{-1} J \left\{ [\hat{P},\psi_2] (1-\Phi^2)(\hat{P})(\hat{P}-z)^{-1} [\hat{P},\psi_1] \right\} J_0 (\hat{P}_T-z)^{-1} J \varphi . \label{halfholomorphic}
\end{eqnarray}
The interest of this decomposition is one hand that $ [\hat{P},\psi_2] \Phi (\hat{P}) $ and $ \Phi (\hat{P}) [\hat{P},\psi_1] $ are bounded ({\it i.e.} have bounded closures) on $ L^2 ({\mathcal M}) $ and on the other hand that the operator $ \{ \cdots \} $ in (\ref{halfholomorphic}) is bounded and holomorphic with respect to $z$ for $ |\mbox{Re}(z)| $ small enough.
Moreover, all these operators have a fast spatial decay in the sense that, for any fixed $ M \geq 1 $, there exist  $ B_1 , B_2 , B (z)  $ bounded on $ L^2 ({\mathcal M}) $ with $ B (z) $ depending boundedly and holomorphically on $z$ for $ \mbox{Re}(z) $ small, such that
$$ [\hat{P},\psi_2] \Phi (\hat{P}) = \scal{r}^{-M} B_2 \scal{r}^{-M},  \qquad \Phi (\hat{P}) [\hat{P},\psi_1] = \scal{r}^{-M} B_1 \scal{r}^{-M} , $$
and
$$  [\hat{P},\psi_2] (1-\Phi^2)(\hat{P})(\hat{P}-z)^{-1} [\hat{P},\psi_1]  = \scal{r}^{-M} B (z) \scal{r}^{-M} . $$
This follows from the compact support of $ [\hat{P},\psi_1] $ and $ [\hat{P},\psi_2] $ and standard arguments (using for instance the Helffer-Sj\"ostrand formula).
Setting
\begin{eqnarray*}
 F (z)  & = & \varphi J_0 (\hat{P}_T -z)^{-1} J \varphi - \varphi J_0 (\hat{P}_T - z)^{-1} \big[ \hat{P}_T,\psi_1 \big] (\hat{P}_T - z)^{-1} J \varphi \nonumber \\
& & - \ \varphi J_0 (\hat{P}_T-z)^{-1} J \scal{r}^{-M} B (z) \scal{r}^{-M}  J_0 (\hat{P}_T-z)^{-1} J \varphi ,
\end{eqnarray*} 
it follows from (\ref{secondline}) and (\ref{partitiondelunitespatiale}) that
\begin{eqnarray}
 (\hat{P}-z)^{-1} &  = & F (z) + (1-\varphi) (\hat{P}-z)^{-1} + \varphi (\hat{P}-z)^{-1} (1-  \varphi )   \nonumber \\
 & & - \varphi J_0 (\hat{P}_T-z)^{-1} J \scal{r}^{-M}B_2 \scal{r}^{-M} (\hat{P}-z)^{-1} \scal{r}^{-M} B_1 \scal{r}^{-M} J_0 (\hat{P}_T-z)^{-1} J \varphi   .
 \nonumber
\end{eqnarray}
We now use this expression to study a spectrally localized version of $ (\hat{P}-z)^{-1} $. We consider
 $\phi_1 (\hat{P}) (\hat{P}-z)^{-1} \phi_2 (\hat{P}) $ where we let $ \phi_1 $ and $ \phi_2 $ take the form
$$ \phi_1 (\lambda) = \Phi (\lambda / \varepsilon_1), \qquad \phi_2 = \Phi (\lambda / \varepsilon_2) . $$
By Lemma \ref{lemmeacitermieux}, we can choose first $ \varepsilon_1 $ small enough such that
$$ || \Phi (\hat{P}/\varepsilon_1) (1-\varphi) \scal{r}^M || \leq 1/4 . $$
Once $ \varepsilon_1 $ is chosen, we can use again Lemma \ref{lemmeacitermieux} to pick $ \varepsilon_2 > 0 $ and $ \nu > 0 $ such that
$$ \big| \big| \scal{r}^{-M} \Phi(\hat{P}/\varepsilon_1) \varphi \scal{r}^M \big| \big| \ \big| \big|  \scal{r}^M (1-\varphi) \Phi (\hat{P}/\varepsilon_2) \big| \big| \leq 1/4 , $$
and
$$ \big| \big| \scal{r}^{-M} \Phi(\hat{P}/\varepsilon_1)\varphi J_0 (\hat{P}_T-z)^{-1} J \scal{r}^{-M}B_2 \big| \big| \ \big| \big| B_1 \scal{r}^{-M} J_0 (\hat{P}_T-z)^{-1} J \varphi \Phi (\hat{P}/\varepsilon_2) (\hat{P}) \scal{r}^{-M} \big| \big| \leq 1/4, $$
for $ |\mbox{Re}(z)| < \nu $ (recall that the first factor in the left-hand side is bounded uniformly with respect to $|\mbox{Re}(z)| \leq 1 $  by (\ref{1ertheoreme})).
Summing up, by choosing
\begin{eqnarray*}
\begin{array}{rclcrcl}
 S_1 & = & \scal{r}^{-M} \Phi (\hat{P}/\varepsilon_1) (1-\varphi) \scal{r}^M, & \qquad & T_1 & = & I , \\
 T_2 & = & \scal{r}^{-M} \Phi(\hat{P}/\varepsilon_1) \varphi \scal{r}^M, & \qquad & S_2 & = & \scal{r}^M (1-\varphi) \Phi (\hat{P}/\varepsilon_2) \scal{r}^{-M} 
 \end{array}
\end{eqnarray*}
and
\begin{eqnarray*}
T (z) & = & \scal{r}^{-M} \Phi(\hat{P}/\varepsilon_1)\varphi J_0 (\hat{P}_T-z)^{-1} J \scal{r}^{-M}B_2 , \\
S(z) & = & B_1 \scal{r}^{-M} J_0 (\hat{P}_T-z)^{-1} J \varphi \Phi (\hat{P}/\varepsilon_2) (\hat{P}) \scal{r}^{-M}  
\end{eqnarray*} 
which satisfy (\ref{smallness}),  we have shown that
\begin{eqnarray*}
 \scal{r}^{-M} (\hat{P}-z)^{-1} \phi_1 (\hat{P}) \phi_2 (\hat{P}) \scal{r}^{-M} & = &  S_1 R_M (z) T_1 + T_2 R_M (z) S_2 + T (z) R_M (z) S (z) \\
 &  & + \ \scal{r}^{-M} \phi_1 (\hat{P}) F (z) \phi_2 (\hat{P}) \scal{r}^{-M} .
\end{eqnarray*}
If we restrict $ z $ to $ |\mbox{Re}(z)| < \varepsilon_0 $ with $ \varepsilon_0 $ small enough, then $ (\hat{P}-z)^{-1} \phi_1 (\hat{P}) \phi_2 (\hat{P}) - (\hat{P}-z)^{-1}  $ becomes holomorphic and bounded in this strip. Therefore, upon adding a bounded holomorphic term to the right-hand side, we may replace the left-hand side by $ \scal{r}^{-M}R_M(z) \scal{r}^{-M} $ and (\ref{neumann}) holds with
$$ F_M (z) = \scal{r}^{-M} \left( (1-\phi_1 \phi_2) (\hat{P}) (\hat{P}-z)^{-1} + \phi_1 (\hat{P}) F (z) \phi_2 (\hat{P}) \right) \scal{r}^{-M} . $$
 The estimates (\ref{borne0}) and (\ref{borne1}) follow directly from  (\ref{1ertheoreme}), (\ref{pourcompacite}) and Lemma \ref{souslemme1}. The estimates (\ref{borne2}) follow from  Lemma \ref{souslemme2}. \finpreuve

\bigskip

 To derive estimates with sharper weights, we will use the following proposition.

\begin{prop} Let $ (w_{\epsilon})_{\epsilon \in (0,1]} $ be the family of weights 
$$ w_{\epsilon} (r) = \scal{r}^{-\mu_1} \scal{\epsilon r}^ {-\mu_2} $$
where $ \mu_1 , \mu_2 \geq 0 $ are fixed real numbers. Fix $ \varphi \in C^{\infty} ({\mathcal M}) $ supported in $ \{r > R \} $, equal to $ 1 $ near infinity,  and let us define
$$ R_{w_{\epsilon}} (z) = w_{\epsilon} (\hat{P}-z)^{-1} w_{\epsilon}, \qquad R_{w_{\epsilon},\varphi}(z) = w_{\epsilon} \varphi (\hat{P}-z)^{-1} \varphi w_{\epsilon} . $$
Then, there exist
\begin{enumerate}
\item{families of bounded operators $ (A_{\epsilon,j} )_{\epsilon \in (0,1]} , (B_{\epsilon,j})_{\epsilon \in (0,1]} , (S_{\epsilon,j})_{\epsilon \in (0,1]}  $ on $ L^2 ({\mathcal M}) $, $ j = 1 , 2 $, such that, for all $ \epsilon \in (0,1] $,
\begin{eqnarray}
 || B_{\epsilon,1} || \ ||S_{\epsilon,1}|| + 
 || B_{\epsilon,2} || \ ||S_{\epsilon,2} || & < & \frac{1}{2} , \label{pourabsorber1} \\
  || A_{\epsilon,1} || + ||A_{\epsilon,1}|| & \leq & C , \label{borneA}   
\end{eqnarray}}
 \item{a real number $  \varepsilon_0  > 0 $ and a bounded holomorphic mapping $ z \mapsto K (z) \in {\mathcal B} (L^2({\mathcal M})) $ defined in a neighborhood of $ \{ |z| \leq \varepsilon_0 \} $,}
\end{enumerate} 
  such that, for all $ z \in \Ca \setminus \Ra $ satisfying $ |z| < \varepsilon_0 $ and all $ \epsilon \in (0,1] $,
\begin{eqnarray} 
R_{w_{\epsilon}} (z) & = & w_{\epsilon} K(z) w_{\epsilon} + A_{\epsilon,2}  R_{w_{\epsilon},\varphi}(z) A_{\epsilon,1}  + B_{\epsilon,1} R_{w_{\epsilon}} (z) S_{\epsilon,1} + S_{\epsilon,2} R_{w_{\epsilon}} (z) B_{\epsilon,2} . \label{reductionresolvante}
\end{eqnarray}
In particular, for all  $ k \geq 0 $, there exists $ C_k $ such that
\begin{eqnarray}
\big| \big| \partial_z^k R_{w_{\epsilon}} (z) \big| \big| \leq C_k \left(1 + \big| \big| \partial_z^k R_{w_{\epsilon},\varphi} (z) \big| \big| \right) , \label{reductionestimation}
\end{eqnarray}
for all $ z \in \Ca \setminus \Ra $ such that $ |z| < \varepsilon_0 $ and all $ \epsilon \in (0,1] $.
\end{prop}

The main interest of this proposition is (\ref{reductionestimation}) which allows to estimate the full resolvent by its cutoff near infinity. This will allow to use (\ref{secondline}) in a convenient way. 

\bigskip

\noindent {\it Proof.} The proof is similar, and simpler, than that of Proposition \ref{propvide}. Let us set $ R (z) = (\hat{P}-z)^{-1} $. For functions $ \phi_1 , \phi_2 \in C_0^{\infty}(\Ra) $ equal to $1$ near zero and to be chosen later, we have  
$$ R(z) \phi_1 (\hat{P}) =  R (z) \varphi \phi_1(\hat{P}) + R (z) (1- \varphi ) \phi_1 (\hat{P}) , $$
where
$$  R (z) \varphi \phi_1(\hat{P}) = (1-\phi_2)(\hat{P}) R (z)  \varphi \phi_1(\hat{P}) + \phi_2 (\hat{P}) R (z)  \varphi \phi_1(\hat{P}) . $$
The second term in the right-hand side of the last formula reads
$$  \phi_2 (\hat{P}) (1- \varphi ) R (z)  \varphi \phi_1(\hat{P}) + \phi_2 (\hat{P}) \varphi R (z)  \varphi \phi_1(\hat{P}) . $$
 Setting for simplicity $ W_{\epsilon} = w_{\epsilon}^{-1} $, and writing $ R (z) = R (z) \phi_1 (\hat{P}) + R (z) (1-\phi_1)(\hat{P}) $ whose second term will be holomorphic with respect to $z$ in a vertical strip around $0$, we find that (\ref{reductionresolvante}) holds with
\begin{eqnarray*}
\begin{array}{rclcrcl}
 B_{\epsilon,1} & = & I, & \qquad & S_{\epsilon,1} & = & W_{\epsilon} (1-\varphi) \phi_1 (\hat{P}) w_{\epsilon}, \\
 B_{\epsilon,2} & =  & W_{\epsilon} \varphi \phi_1(\hat{P}) w_{\epsilon}  , & \qquad & S_{\epsilon,2} &  = & w_{\epsilon} \phi_2 (\hat{P}) (1-\varphi) W_{\epsilon} , \\
 A_{\epsilon,1} & = & W_{\epsilon} \phi_1(\hat{P}) w_{\epsilon} , & \qquad & A_{\epsilon,2} & = & w_{\epsilon} \phi_2 (\hat{P}) W_{\epsilon} ,
 \end{array}
\end{eqnarray*}
and
$$ K (z) =  R (z) (1-\phi_1)(\hat{P}) + (1-\phi_2)(\hat{P}) R (z) \varphi \phi_1(\hat{P})  . $$
We choose $ \phi_1 , \phi_2 $ successively as follows.  By using the uniform boundedness of $ W_{\epsilon} \scal{r}^{-\mu_1 - \mu_2} $ and of $ w_{\epsilon} $, 
we can choose first $ \phi_1 (\lambda) = \Phi (\lambda / \varepsilon_1) $
 such that, by using (\ref{compactfacile}) with $ M = \mu_1 + \mu_2 $, we have
$||B_{\epsilon,1}|| \ || S_{\epsilon,1} || < 1/4 $ (recall that $ 1 - \varphi $ is a compactly support spatial cutoff). Then $ B_{\epsilon,2} $, which depends on $ \phi_1 $,  is bounded on $ L^2 ({\mathcal M}) $, uniformly in $ \epsilon $, by routine arguments. We choose then $ \phi_2 (\lambda) = \Phi (\lambda/\varepsilon_2) $ such that the norm of  $ S_{\epsilon,2} $ is small enough to guarantee that $ || B_{\epsilon,2} || ||S_{\epsilon,2} || < 1/4 $. We therefore get (\ref{pourabsorber1}). The operators $ A_{\epsilon,1} $ and $A_{\epsilon,2} $ are uniformly bounded on $ L^2 ({\mathcal M}) $ similarly to $ B_{\epsilon,2} $, which yields (\ref{borneA}). We then fix $ \varepsilon_0 $ small enough such that both $ \phi_1 $ and $ \phi_2 $ are equal to $1$ near $ [-\varepsilon_0,\varepsilon_0] $. This implies that $ K (z) $ is holomorphic and bounded near the strip $ \{|\mbox{Re}(z)| \leq \varepsilon_0 \} $. Finally, (\ref{reductionestimation}) is an easy consequence of (\ref{pourabsorber1}), (\ref{borneA}) and (\ref{reductionresolvante}) and the fact that $ || \partial_z^k K (z) || $ is bounded on $ \{| \mbox{Re}(z)| \leq \varepsilon_0 \} $. 
\finpreuve

\bigskip

We now apply the results of this section to prove Theorems \ref{sharpweightresolvent},  
\ref{pseudorescaledestimates} and \ref{carreresolvante}. The three proofs follow the same strategy, but we present them separately for pedagogical reasons: we consider first the simplest case which is Theorem \ref{sharpweightresolvent} and then explain how to modify the arguments for the other two.

\bigskip

\noindent {\bf Proof of Theorem \ref{sharpweightresolvent}.} We prove first a uniform estimate in $z$, but with a rough weight. By  Proposition \ref{propvide}, we have
$$ \big| \big| \scal{r}^{-M} (\hat{P} - z)^{-1} \scal{r}^{-M} \big| \big| \leq C + \frac{3}{4} \big| \big| \scal{r}^{-M} (\hat{P} - z)^{-1} \scal{r}^{-M} \big| \big| , $$
for all $z \in \Ca \setminus \Ra$ such that $ |\mbox{Re}(z)| < \varepsilon_0 $. This implies that 
\begin{eqnarray} 
 \big| \big| \scal{r}^{-M} (\hat{P} - z)^{-1} \scal{r}^{-M} \big| \big| \leq 4 C . \label{pouroptimalite}
\end{eqnarray}
We now use (\ref{secondline}) to replace $ \scal{r}^{-M} $ by the optimal weight $ \scal{r}^{-1} $ as follows. We start by observing that (\ref{pouroptimalite}) and the compact support of $ [\hat{P},\psi_1] $ and $ [\hat{P},\psi_2] $ in (\ref{secondline}) imply that
\begin{eqnarray} 
 \big| \big| \scal{r} [\hat{P},\psi_2] (\hat{P} - z)^{-1} [\hat{P},\psi_1] \scal{r} \big| \big| \leq  C^{\prime} . \label{pouroptimalite2}
\end{eqnarray}
The weight $ \scal{r} $ could be replace by any power of $ \scal{r} $ but (\ref{pouroptimalite2}) will be sufficient.
Here the unboundeness of the operators $ [\hat{P},\psi_1] $ and $ [\hat{P},\psi_2] $ can easily be overcome for instance by writing
$$ [\hat{P},\psi_2] (\hat{P} - z)^{-1} [\hat{P},\psi_1] = [\hat{P},\psi_2] \Phi (\hat{P}) (\hat{P} - z)^{-1} \Phi (\hat{P}) [\hat{P},\psi_1] + [\hat{P},\psi_2] (\hat{P} - z)^{-1}(1-\Phi^2(\hat{P})) [\hat{P},\psi_1] , $$
with $ \Phi \in C_0^{\infty}(\Ra) $, $ \Phi \equiv 1 $ near $ [-\varepsilon_0, \varepsilon_0] $. Both $ \scal{r} [\hat{P},\psi_2] \Phi (\hat{P}) \scal{r}^{M} $ and  $ \scal{r}^M \Phi (\hat{P}) [\hat{P},\psi_1] \scal{r} $ are bounded on $ L^2 ({\mathcal M}) $ and the second term is holomorphic and bounded near the strip $ \{ |\mbox{Re}(z)| \leq \varepsilon_0 \} $, so (\ref{pouroptimalite2}) follows clearly from (\ref{pouroptimalite}). 
Then, by composing with $ \scal{r}^{-1} $ to the left and to the right of both sides of   (\ref{secondline}), we obtain
\begin{eqnarray*}
  \big| \big| \scal{r}^{-1} \varphi (\hat{P} - z)^{-1} \varphi \scal{r}^{-1} \big| \big| & \lesssim &   \big| \big| \scal{r}^{-1}  (\hat{P}_T - z)^{-1} \scal{r}^{-1}  \big| \big| \\
& & \  + \  \big| \big| \scal{r}^{-1}  (\hat{P}_T - z)^{-1} [\hat{P}_T,\psi_1] (\hat{P}_T-z)^{-1} \scal{r}^{-1}  \big| \big| \\
  &  & \ + \  \big| \big| \scal{r}^{-1}  (\hat{P}_T - z)^{-1} \scal{r}^{-1}  \big| \big|^2  , 
\end{eqnarray*}  
where last term (squared) is the contribution of the last term of (\ref{secondline}) combined with (\ref{pouroptimalite2}).
The second term in the right-hand side can be estimated thanks to Lemma \ref{souslemme1} with $ w = \scal{r}^{-1} $ and $k=0$
so, using (\ref{1ertheoreme}), we conclude that $ \big| \big| \scal{r}^{-1} \varphi (\hat{P} - z)^{-1} \varphi \scal{r}^{-1} \big| \big| $ is uniformly bounded for $ |\mbox{Re}(z)| $ small. We can then remove the cutoff $ \varphi $ by using (\ref{reductionestimation}) (with $ \mu_2 =0 $ and $ \mu_1 = 1 $) and this completes the proof.
\finpreuve

\bigskip

Since its proof is less technical, we prove Theorem \ref{carreresolvante} before Theorem \ref{pseudorescaledestimates}.

\bigskip 
 
 \noindent {\bf Proof of Theorem \ref{carreresolvante}.} By differentiating (\ref{neumann}) with respect to $z$, we obtain 
$$ \big| \big| \partial_z R_M (z) \big| \big| \leq \frac{3}{4} \big| \big| \partial_z R_M (z) \big| \big| + \big| \big| \partial_z F_M (z) \big| \big|
+ C \big| \big| \partial_z T_M (z) \big| \big| + C \big| \big| \partial_z S_M (z) \big| \big| , $$
using (\ref{pouroptimalite}) to control the contribution of the non differentiated $ R_M (z)  $ in the right-hand side. Using (\ref{borne1}), this implies that
$$   \big| \big| \scal{r}^{-M} (\hat{P}-z)^{-2} \scal{r}^{-M} \big| \big| = \big| \big| \partial_z R_M (z) \big| \big|   \lesssim |\mbox{Re}(z)|^{s-1}, $$
when $ |\mbox{Re}(z)|>0 $ is small enough. It remains to replace the weight $ \scal{r}^{-M} $ by $ \scal{r}^{-2-s} $. The first observation is that, for the very same reason we got (\ref{pouroptimalite2}), the above estimate implies that
\begin{eqnarray}                           
\big| \big| \scal{r}^{M} [\hat{P},\psi_2] (\hat{P} - z)^{-2} [\hat{P},\psi_1] \scal{r}^{M} \big| \big| \lesssim |\mbox{Re}(z)|^{s-1} . 
 \label{aswellas}     
 \end{eqnarray}                      
Using (\ref{1ertheoreme}), (\ref{pourcompacite}) and Lemma \ref{souslemme1} with $ w = \scal{r}^{-2-s} $, we also have
\begin{eqnarray}
 \big| \big| \partial_z \big(  \scal{r}^{-2-s}  (\hat{P}_T - z)^{-1} [\hat{P}_T,\psi_1] (\hat{P}_T-z)^{-1} \scal{r}^{-2-s} \big) \big| \big|  \lesssim |\mbox{Re}(z)|^{s-1} .
 \label{similarlyzerobis}
\end{eqnarray}
Therefore, it follows from (\ref{secondline}) that
$$   \big| \big| \scal{r}^{-2-s} \varphi (\hat{P}-z)^{-2} \varphi \scal{r}^{-2-s} \big| \big| \lesssim |\mbox{Re}(z)|^{s-1} , $$
using again (\ref{pourcompacite}) for the first term as well as (\ref{pouroptimalite}) and (\ref{aswellas}) to handle the contribution of the last term. We can then remove the factors $ \varphi $ by using (\ref{reductionestimation}) with $ \mu_2 = 0 $ and $ \mu_1 = 2 +s $. The result follows. \finpreuve

\bigskip

\noindent {\bf Proof of Theorem \ref{pseudorescaledestimates}.} Using Proposition \ref{propvide}, it is not hard to check by (finite) induction on $k$ that
\begin{eqnarray}
 || \partial_z^k R_M (z) || \lesssim |\mbox{Re}(z)|^{-k} , \label{atraduire}
\end{eqnarray}
for $ 0 < \mbox{Re}(z) < \varepsilon_0 $. Using that
$$ (\epsilon^{-2} \hat{P} - Z)^{-1-k} = \frac{\epsilon^{2}}{k!} \partial_Z^k (\hat{P} - \epsilon^2 Z)^{-1} = \left. \epsilon^{2(k+1)} \frac{\partial_z^k}{k!} (\hat{P}-z)^{-1} \right|_{z = \epsilon^2 Z} $$
we obtain from (\ref{atraduire}) the a priori estimates
\begin{eqnarray}
 || \scal{r}^{-M} (\epsilon^{-2}\hat{P}-Z)^{-1-k} \scal{r}^{-M} || \lesssim \epsilon^{2(k+1)} |\mbox{Re}(\epsilon^2 Z)|^{-k} \lesssim \epsilon^2 . \label{238}
\end{eqnarray}
On the other hand, it is not hard to check that (\ref{estimeesderniertheoreme}) is equivalent to
\begin{eqnarray}
 || \scal{\epsilon r}^{-1-k} \partial_Z^k (\hat{P} - \epsilon^2 Z)^{-1} \scal{\epsilon r}^{-1-k} || \lesssim \epsilon^{-2}  ,   \label{citationfinale}
\end{eqnarray}
with $ k = N - 1 $. We will   show that (\ref{citationfinale}) holds for all $k$ between $ 0 $ and $ N-1$. By (\ref{238}), we already know that
\begin{eqnarray}
 || \scal{ r}^{-M} \partial_Z^k (\hat{P} - \epsilon^2 Z)^{-1} \scal{ r}^{-M} || \lesssim 1 ,\label{238bis}
\end{eqnarray}
which are better estimates, but with the much stronger $ \epsilon $ independent weight $ \scal{r}^{-M} $. We shall 
 combine (\ref{238bis})  with a priori estimates on $ (\epsilon^{-2}\hat{P}_T-Z)^{-1} $ to derive (\ref{estimeesderniertheoreme}). Our first observation is that, 
by (\ref{rescalezero}), we have slightly more precise estimates for $ \hat{P}_T $, namely
\begin{eqnarray}
 || \scal{\epsilon r}^{-k-1} \partial_Z^k (\hat{P}_T - \epsilon^2 Z)^{-1} \scal{\epsilon r}^{-k-1} || & \lesssim & \epsilon^{-2} , \label{aciteraussi} \\
 || \scal{\epsilon r}^{-k-1} \partial_Z^k (\hat{P}_T - \epsilon^2 Z)^{-1} \scal{\epsilon r}^{-k} \scal{r}^{-1} || & \lesssim & \epsilon^{-1} , \label{intermediairerescale} 
\end{eqnarray}
This follows from (\ref{rescalezero}) since we are allowed to use one power of the homogeneous weight $ r^{-1} $ at most on each side. 
For instance, (\ref{intermediairerescale}) is obtained by
\begin{eqnarray*}
 || \scal{\epsilon r}^{-k-1} \partial_Z^k (\hat{P}_T - \epsilon^2 Z)^{-1} \scal{\epsilon r}^{-k} \scal{r}^{-1} || & \lesssim & 
 \epsilon^{2k} || \scal{\epsilon r}^{-k-1}  (\hat{P}_T - \epsilon^2 Z)^{-1-k} \scal{\epsilon r}^{-k} \scal{r}^{-1} || \\
 & \lesssim & \epsilon^{-2}  || \scal{\epsilon r}^{-k-1}  (\epsilon^{-2}\hat{P}_T -  Z)^{-1-k} \scal{\epsilon r}^{-k} \scal{r}^{-1} || \\
 & \lesssim & \epsilon^{-1}  || \scal{\epsilon r}^{-k-1}  (\epsilon^{-2}\hat{P}_T -  Z)^{-1-k} \scal{\epsilon r}^{-k} (\epsilon r)^{-1} || \\
 & \lesssim & \epsilon^{-1} ,
\end{eqnarray*} 
where we  pass from the second to the third line by replacing $ \scal{r}^{-1} $ by $ r^{-1} $, and get the final estimate by using (\ref{rescalezero}). 
Using (\ref{aciteraussi}), (\ref{intermediairerescale}) and Lemma \ref{souslemme1},
it follows that
\begin{eqnarray}
  \big| \big| \partial_Z^k \big( \scal{\epsilon r}^{-k-1} (\hat{P}_T- \epsilon^2 Z)^{-1} [\hat{P}_T,\psi] (\hat{P}_T- \epsilon^2Z)^{-1} \scal{\epsilon r}^{-k-1} \big)  \big| \big| \lesssim \epsilon^{-2} . \label{aciteraussibis} 
\end{eqnarray}
Also, using (\ref{238bis}) and  (\ref{intermediairerescale}), the contribution of the last term of (\ref{secondline}) is
\begin{eqnarray}
 \big| \big| \partial_Z^k \big( \scal{\epsilon r}^{-k-1} J_0 (\hat{P}_T- \epsilon^2 Z)^{-1} J [\hat{P},\psi_2] (\hat{P}- \epsilon^2 Z)^{-1} [\hat{P},\psi_1] J_0 (\hat{P}_T-\epsilon^2 Z)^{-1} J  
  \scal{\epsilon r}^{-k-1} \big)  \big| \big| \lesssim \epsilon^{-2} .                        \nonumber            
  \end{eqnarray}
 Therefore, this last estimate together with (\ref{aciteraussi}), (\ref{aciteraussibis}) and  (\ref{secondline}) 
  imply  that for all $ k $ between $ 0 $ and $ N - 1 $, 
 $$ \big| \big| \partial_Z^k \big( \scal{\epsilon r}^{-k-1} \varphi (\hat{P}- \epsilon^2 Z)^{-1}  \varphi \scal{\epsilon r}^{-k-1} \big)  \big| \big| \lesssim \epsilon^{-2} . $$
 We can then drop the cutoff $ \varphi $ by using (\ref{reductionestimation}). We have thus proved (\ref{citationfinale}) which, in the special case $ k = N-1 $, yields (\ref{estimeesderniertheoreme}).  \finpreuve

\section{Mourre theory on a cone} \label{Mourreonacone}
\setcounter{equation}{0}

In this  section, we develop a Mourre theory for elliptic operators in divergence form on an exact cone, which will be crucial to prove the estimates (\ref{1ertheoreme}), (\ref{rescalezero}) and (\ref{pourcompacite}).
\subsection{Operators in divergence form}
We start by introducing a class of tensors which will be convenient to handle operators in divergence form on the cone $ {\mathcal M}_0 = (0,+\infty) \times {\mathcal S} $ introduced in (\ref{Mzero}). Given a Riemannian metric $ g $ on  $ {\mathcal M}_0 $, we denote the associated co-metric ({\it i.e.} the inner product on the fibers of $ T^* {\mathcal M}_0 $) by $ g^* $.  Then there exists a unique tensor $ T^g \in C^{\infty} \big( {\mathcal M}_0 , \mbox{Hom}(T^* {\mathcal M}_0 , T{\mathcal M}_0) \big)   $, {\it i.e.} a section of the vector bundle $ \mbox{Hom}(T^* {\mathcal M}_0 , T{\mathcal M}_0) \approx (T^* {\mathcal M}_0)^* \otimes T {\mathcal M}_0 $, such that for all $1 $-forms $ \xi, \eta $ on $ {\mathcal M}_0 $, 
\begin{eqnarray} 
 \eta \cdot T^g \xi  = g^* (\xi , \eta) , \label{generaldefinition}
\end{eqnarray}
where $ \cdot $ is the intrinsinc duality between  a $1$-form and a vector. In usual terms, $ T^g$ raises indices. 
It is automatically symmetric in the sense that we have 
\begin{eqnarray}
  \xi \cdot T^g \eta = \eta \cdot  T^g  \xi   , \label{tenseursymmetrique}
\end{eqnarray} 
for all $ \xi , \eta $,
 and it allows to write the Laplace-Beltrami operator $ \Delta_g $ as
$$ \Delta_g u =  \mbox{div}_g \big( T^g d u \big) , $$
for all smooth functions $ u  $.

Using the isomorphisms $ T {\mathcal M}_0 \approx T \Ra^+ \times T {\mathcal S} $ and $ T^* {\mathcal M}_0 \approx T^* \Ra^+ \times T^* {\mathcal S} $, we can write any tensor $ T  \in C^{\infty} \big( {\mathcal M}_0 , \mbox{Hom}(T^* {\mathcal M}_0 , T{\mathcal M}_0) \big) $  in matrix form as
$$ T = \left( \begin{matrix} T_{11} & T_{12}  \\ T_{21} & T_{22}  \end{matrix} \right) $$
with $ T_{11} \in C^{\infty} (\Ra^+ \times {\mathcal S}) $, $ T_{22} \in C^{\infty} \big(\Ra^+ , C^{\infty} \big( {\mathcal S} , \mbox{Hom}(T^* {\mathcal S} , T {\mathcal S} ) \big) \big) $
and
$$ T_{12} \in C^{\infty} \big(\Ra^+ , C^{\infty} \big( {\mathcal S} , \mbox{Hom}(T^* {\mathcal S} , \Ra) \big) \big), \qquad T_{21} \in 
C^{\infty} \big(\Ra^+ , C^{\infty} \big( {\mathcal S} , \mbox{Hom}(\Ra , T {\mathcal S}) \big) \big) . $$
For all $ \omega \in {\mathcal S} $, any element of $ u_{\omega}  \in \mbox{Hom}(\Ra , T_{\omega} {\mathcal S}) $ has an intrinsic adjoint (or transpose) denoted by  $ u^{\dag}_{\omega} \in \mbox{Hom}( T^*_{\omega} {\mathcal S},\Ra) $ and defined by
$$ u^{\dag}_{\omega} (\xi) = \xi \cdot u_{\omega} , \qquad \xi \in T_{\omega}^* {\mathcal S} . $$
If $ u \in C^{\infty} \big( {\mathcal S} , \mbox{Hom}(\Ra , T {\mathcal S}) \big) $ is a section, we define the section $ u^{\dag} \in C^{\infty} \big( {\mathcal S} ,  \mbox{Hom} (T^* {\mathcal S} , \Ra) ) \big) $ in the obvious way ($ u^{\dag}(\omega) = u (\omega)^{\dag} $).
It is then easy to check  the following characterization of symmetric tensors, in the sense of (\ref{tenseursymmetrique}),
$$T \ \mbox{is symmetric} \ \ \Longleftrightarrow \ \ \mbox{for each} \ r \in \Ra^+, \ \ T_{12}(r) = T_{21}(r)^{\dag} \ \ \mbox{and} \ T_{22}(r) \ \mbox{is symmetric} $$
where, for each $r$, the tensors in the right-hand side belong respectively to $ C^{\infty} \big( {\mathcal S} , \mbox{Hom} (T^* {\mathcal S} , \Ra) \big) $ and $ C^{\infty} \big( {\mathcal S} , \mbox{Hom} (T^* {\mathcal S} , T {\mathcal S}) \big) $.

If we consider the  conical metric
\begin{eqnarray}
 g_0 = d r^2 + r^2 h_0 , \label{exactconicalmetric}
\end{eqnarray} 
   on $  {\mathcal M}_0 $ and use the above formalism, we then have
$$ T^{g_0}  = \left( \begin{matrix} 1 & 0 \\ 0 & r^{-1} I_{T {\mathcal S}} \end{matrix} \right) 
\left( \begin{matrix} 1 & 0 \\ 0 & T^{h_{0}} \end{matrix} \right) \left( \begin{matrix} 1 & 0 \\ 0 & r^{-1} I_{T^* {\mathcal S}} \end{matrix} \right) . $$
We now introduce a class of perturbations of this tensor.  
For any $ V \in C^{\infty}\big( {\mathcal S} , \mbox{Hom} (\Ra, T {\mathcal S}) \big) $, which can be identified with a vector field on $ {\mathcal S} $, we set
\begin{eqnarray}
 || V ||_{L^{\infty}({\mathcal S})} := \sup_{\omega \in {\mathcal S}} |V_{\omega}|_{h_{0,\omega}} , \label{Linfini1}
\end{eqnarray} 
$ | \cdot |_{h_{0,\omega}} $ denoting the norm on $ T_{\omega} {\mathcal S} $ associated to $ h_0 $. If $ W \in C^{\infty} \big( {\mathcal S} , \mbox{Hom} (T^* {\mathcal S}, T {\mathcal S} ) \big) $, we set 
\begin{eqnarray}
 || W_{\omega} ||_{h_0} = \sup  \frac{|\eta \cdot W_{\omega}(\xi) |}{ | \xi|_{h^*_{0,\omega}}|\eta|_{h^*_{0,\omega}} } , \label{Linfini2}
\end{eqnarray} 
  the sup being taken over all $ \xi , \eta \in T_{\omega}^* {\mathcal S} \setminus \{ 0 \} $, and $ | \cdot |_{h_{0,\omega}^*} $ being the norm on $ T^*_{\omega} {\mathcal S} $. Then, we set
$$ || W ||_{L^{\infty}({\mathcal S})} := \sup_{\omega \in {\mathcal S}} || W_{\omega} ||_{h_0} . $$

\begin{defi}[admissible perturbation] \label{aumoinspourlanorme} Fix an integer $ N \geq 0 $. A symmetric tensor $ K \in C^{\infty} \big( {\mathcal M}_0, \emph{Hom}(T^* {\mathcal M}_0 , T{\mathcal M}_0) \big) $ of the form 
$$ K =  
\left( \begin{matrix} K_{11}(r) & K_{21}(r)^{\dag} \\ K_{21}(r)  & K_{22}(r)  \end{matrix} \right)  $$
is an $N$-admissible perturbation if for all $ k \leq N  + 1 $,
$$ ||| K |||_{k} := \sup_{r > 0} \big| \big| (r \partial_r)^k  K_{11}(r) \big| \big|_{L^{\infty}({\mathcal S})} + \big| \big| (r\partial_r)^k K_{21}(r) \big| \big|_{L^{\infty}({\mathcal S})} + \big| \big| (r\partial_r)^k K_{22} (r) \big| \big|_{L^{\infty}({\mathcal S})}  < \infty . $$ 
\end{defi}

We note that although we assume $ K $ to be smooth, we only require a control on the derivatives with respect to $r$. The interest of this class is to behave nicely under rescaling in $r$. For future purposes we record here the notation, 
\begin{eqnarray}
 K^t (r) := K (e^t r ) , \label{Kreechelonne}
\end{eqnarray}
defined for $t \in \Ra$ and any admissible perturbation $K$.

We now consider differential operators associated to such perturbations. In the sequel, we will denote the Riemannian measure associated to $ g_0 $ by,
$$ d \mu = r^{n-1} d r d {\rm vol}_{h_0}  . $$
For a given $ N $-admissible perturbation $ K $, we can consider the sesquilinear form
\begin{eqnarray}
 Q_{K} (u,v) =  \int_{\Ra^+ \times {\mathcal S}} \left( \begin{matrix} \partial_r \overline{u} \\  d_{\mathcal S} \overline{u}/r \end{matrix} \right) \cdot \left\{ \left( \begin{matrix} 1 & 0 \\ 0 & T^{h_0} \end{matrix} \right) + \left( \begin{matrix} K_{11} & K_{21}^{\dag} \\ K_{21}  & K_{22} \end{matrix} \right) \right\} \left( \begin{matrix} \partial_r v \\  d_{\mathcal S} v /r \end{matrix} \right) d \mu , \label{definitionformeK}
\end{eqnarray}
first for $ u , v \in  C_0^{\infty}({\mathcal M}_0) $. For simplicity, everywhere in the sequel we set
$$ C_0^{\infty} = C_0^{\infty}({\mathcal M}_0) . $$ 
If we let $ (.,.)_{L^2} $ be the inner product of $ L^2 ({\mathcal M}_0, d\mu) $, we see by integration by part that
\begin{eqnarray}
 Q_{K}(u,v) =  (P_K u , v)_{L^2} , \nonumber
\end{eqnarray}
with
\begin{eqnarray}
 P_K v = - \Delta_{g_0} v - \mbox{div}_{g_0} \big( K^{\rm sc} d v \big) , \label{definitionPK}
\end{eqnarray}
where 
$$ K^{\rm sc} = \left( \begin{matrix} 1 & 0 \\ 0 & r^{-1} I_{T {\mathcal S}} \end{matrix} \right) 
\left( \begin{matrix} K_{11}(r) & K_{21} (r)^{\dag} \\ K_{21} (r) & K_{22} (r) \end{matrix} \right) \left( \begin{matrix} 1 & 0 \\ 0 & r^{-1} I_{T^* {\mathcal S}} \end{matrix} \right) . $$
Notice that, by the symmetry assumption on $K$, the operator $ P_{K} $ is symmetric with respect to $ d \mu $. 
%
All this will allow to define closed realizations of the differential operators $ P_K $ by means of sesquilinear forms (see Subsection \ref{soussectionMourre}). Before doing so, we need to introduce the relevant Sobolev norms, as well as useful intermediate results in the next subsection. 

\subsection{Sobolev spaces and dilations}
Let us consider the operators
$$ D_r = i^{-1} \partial_r, \qquad  |D_{\mathcal S}| = (- \Delta_{\mathcal S})^{1/2}  $$
that preserve $ C_0^{\infty}  $. We let $ e^{itA} $ be the unitary group defined on $ L^2   $ by
\begin{eqnarray}
 \big( e^{itA} u \big) (r,\omega) = e^{t n / 2} u (e^t r , \omega) , \label{dilations}
\end{eqnarray}
whose generator is the differential operator $ A = \frac{n}{2i} - i r \partial_r $. Note that $ e^{itA} $ preserves $ C_0^{\infty}  $.

\begin{prop} On $ C_0^{\infty}  $, the following identities hold
\begin{eqnarray}
e^{-itA} D_r e^{itA} & = & e^{t} D_r , \label{A6} \\
e^{-itA} r^s e^{itA} & = & e^{-ts}r^s , \qquad s \in \Ra . \label{A7} 
\end{eqnarray}
If $ K $ is an $N$-admissible perturbation and $ K_t $ is given by (\ref{Kreechelonne}), then
\begin{eqnarray}
Q_{K} (e^{-itA} u , v)  & = &  e^{-2t} Q_{K^t} (u , e^{itA} v)   , \label{commutateurparderivee}
\end{eqnarray}
for all $u,v \in C_0^{\infty} $. In particular, if $ N \geq 1 $,
\begin{eqnarray}
i \big( Q_K ( u ,  A v)  - Q_K ( A u ,  v) \big)  =  \big( u ,  i [P_K,A]  v \big)_{L^2} =  2 Q_{K_1} (u, v) \label{commutateurformel0}
\end{eqnarray}
with
$$ K_1 = \left( 1 - \frac{r \partial_r}{2} \right) K . $$
\end{prop}

\noindent {\it Proof.} The formulas (\ref{A6}), (\ref{A7}), (\ref{commutateurparderivee}) are routine and (\ref{commutateurformel0}) follows from 
$$ 2 Q_K (u,v) - \frac{d}{dt} Q_{K^t} (u,v)_{|t=0} = i \big( Q_K (u,Av) - Q_K (Au,v) \big) , $$
by differentiating (\ref{commutateurparderivee}) in $t$ and evaluating it at $ t =0 $. \finpreuve

\bigskip

We next define the norm
$$ || u ||_{H_0^1} = \left( || u ||_{L^2}^2 + || D_r u ||_{L^2}^2 + \big| \big| r^{-1} |D_{\mathcal S}| u  \big|\big|_{L^2}^2 \right)^{1/2} , $$
first on $ C_0^{\infty} $ and then  on $ H_0^1 := H_0^1 ( {\mathcal M}_0) $ defined as
\begin{eqnarray}
 H_0^1 = \mbox{closure of} \ C_0^{\infty} \ \mbox{for the norm} \ || \cdot ||_{H_0^1} . \label{defH01}
\end{eqnarray}
The operators $ D_r $ and $ r^{-1} |D_{\mathcal S}| $ have unique continuous extensions as linear maps from $ H_0^1 $ to $ L^2 $.
It is also convenient to introduce the homogeneous Sobolev norm
$$ || u ||_{\dot{H}_0^1} := \big( || D_r u ||_{L^2}^2 + \big| \big| r^{-1} |D_{\mathcal S}| u \big| \big|_{L^2}^2 \big)^{1/2} , $$
which we shall  consider only on $ H_0^1 $, so we do not need to introduce the corresponding space $ \dot{H}_0^1 $.
We finally set
$$ H^{-1} = \mbox{topological dual space to} \ H_0^1 . $$
We denote the antilinear duality between $ u \in  H_0^1 $ and $ E \in H^{-1} $ by $ ( E , u ) $, with the convention that it is linear in $u$ and conjugate linear in $ E $. In other words, if $ \scal{.,.} $ is the bilinear pairing between $ H_0^1 $ and its dual, we have set
\begin{eqnarray}
 (E,u) := \overline{\scal{E,\bar{u}}} . \label{symetriser0}
\end{eqnarray} 
To make this definition more symmetric, we also set
\begin{eqnarray}
 (u,E) := \overline{\left( E , \overline{u} \right)}, \label{symetriser}
\end{eqnarray} 
  for all $ E \in H^{-1} $ and $ u \in H_0^1 $.  We have the following useful and elementary result which we  record at least for notational purpose.

\begin{prop} \label{injection} For all $ f \in L^2 $, there exists  a unique $ E_f \in H^{-1} $ such that, for all $ u \in H_0^{1} $
$$ (f,u)_{L^2} = (E_f,u ) . $$
The map $ f \mapsto E_f $ is linear, continuous and injective, thus realizes an embedding from $ L^2 $ into $ H^{-1} $. We denote it by $ \bar{I} $. Moreover, $ L^2 $ ({\it i.e.} $\bar{I}L^2$) is dense in $ H^{-1} $.
\end{prop}

We omit the proof which is standard.
The interest of this proposition is to be able to consider $ L^2 $ as a (dense) subspace of $ H^{-1} $. We shall use this convenient identification everywhere in the sequel. For instance, if $ E $ belongs to $ L^2 $, (\ref{symetriser0}) and (\ref{symetriser}) correspond to  $ L^2 $ inner products. 


\bigskip

We next summarize several useful properties on $ H_0^1 $ and $ H^{-1} $ related to the group (\ref{dilations}). We will be in particular  interested
in the properties of the resolvent of $A$, as an operator on $ L^2 $, when restricted to $ H_0^1 $. We recall that, if $ \alpha > 0 $ is a real parameter and $ \zeta \in \Ca \setminus \Ra $, we have
\begin{eqnarray}
 (\alpha A - \zeta)^{-1} = \frac{1}{i} \int_0^{\pm \infty} e^{-it\zeta} e^{it \alpha A} dt, \qquad \mp \mbox{Im}(\zeta) > 0 . \label{formulaA}
\end{eqnarray}

\begin{prop} \label{differentiabiliteH01} \begin{enumerate} 
\item{$ H_0^1 $ is stable by multiplication by smooth functions of $r$ which are bounded together with their derivatives.}
\item{$ H_0^1 $ is stable by $ e^{itA} $, $ e^{itA} $ is strongly continuous on $ H_0^1 $ and 
$$ || e^{itA}  u||_{H_0^1} \lesssim (1 +e^{t}) || u ||_{H_0^1} , \qquad t \in \Ra, \ u \in H_0^1 .$$ Furthermore, if $  0 < \alpha < |\emph{Im}(\zeta)| $, $ H_0^1 $ is stable by $ (\alpha A-\zeta)^{-1} $. }
\item{The group $ e^{itA} $ extends from $ L^2 $ to an $ H^{-1} \rightarrow H^{-1} $ strongly continuous group. Its adjoint is $ e^{-itA} $ (acting on $ H_0^1 $). Furthermore, if $  0 < \alpha < |\emph{Im}(\zeta)| $, $ (\alpha A-\zeta)^{-1} $ extends from $ L^2 $ to a bounded $ H^{-1} \rightarrow H^{-1} $ operator, whose adjoint (acting on $ H_0^1 $) is $ (\alpha A-\bar{\zeta})^{-1} $.}
\item{Fix $ 0 < \alpha < |\emph{Im}(\zeta)|$. For all $ t \in \Ra$, 
$$ e^{itA} (\alpha A-\zeta)^{-1} = (\alpha A-\zeta)^{-1} + \frac{i}{\alpha} \int_0^t e^{isA} \left( I + \zeta (\alpha A-\zeta)^{-1} \right) ds $$
as an equality between operators on $ H_0^1 $ (resp. $ H^{-1} $). Here the integral converges in the strong sense. In particular $ e^{itA} (\alpha A-\zeta)^{-1} $ is strongly differentiable with respect to $t$ on $ H_0^1 $ and $ H^{-1} $.} 
\item{Fix $ \alpha \in (0,1) $ and an integer $ N \geq 1 $. We have the interpolation estimate
$$ || A (\alpha A+i)^{-N} u ||_{H_0^1} \leq C || (\alpha A+i)^{-N} u ||_{H_0^1}^{1-\frac{1}{N}} || u ||_{H_0^1}^{\frac{1}{N}} , $$
for all $ u \in H_0^1 $.}  
\end{enumerate}
\end{prop}

\bigskip

\noindent {\it Proof.} The item 1 is straightforward by density of $ C_0^{\infty} $. The estimate of the item 2 holds true on $ C_0^{\infty} $ by (\ref{A6}) and (\ref{A7}) (with $ s=-1 $) hence on $ H_0^1 $ by density. Checking the strong continuity is a routine. The boundedness of $ (\alpha A-\zeta)^{-1} $ on $ H_0^1 $ is then a consequence of (\ref{formulaA}).
  The item 3 follows from Proposition \ref{injection}, the item 2 of the present proposition and  the formula (\ref{formulaA}) combined with routine duality arguments.
  The identity of the item 4 holds clearly on $ L^2 $ since $ e^{itA} (\alpha A-\zeta)^{-1} $ is strongly $ C^1 $ in $t$ (note that $ i \alpha^{-1}e^{isA} \left( I + \zeta (\alpha A-\zeta)^{-1} \right) = e^{isA} i A (\alpha A-\zeta)^{-1}$). That the integral converges in the strong sense on $ H_0^1 $ (resp. $ H^{-1} $) follows from the item 2 (resp. item 3).
To prove the item 5, we recall first that
\begin{eqnarray}
 || \scal{A} (\alpha A+i)^{-N}  f ||_{L^2} \leq C || (\alpha A+i)^{-N} f ||_{L^2}^{1-\frac{1}{N}} || f ||_{L^2}^{\frac{1}{N}} , \label{interpolationA}
\end{eqnarray} 
using the spectral theorem and the Hadamard three lines theorem.
In particular, we have
\begin{eqnarray*}
 || A (\alpha A+i)^{-N}  u ||_{L^2}  & \leq & C || (\alpha A+i)^{-N} u ||_{L^2}^{1-\frac{1}{N}} || u ||_{L^2}^{\frac{1}{N}} \\
   & \leq & C || (\alpha A+i)^{-N} u ||_{H_0^1}^{1-\frac{1}{N}} || u ||_{H_0^1}^{\frac{1}{N}} .
\end{eqnarray*}
It remains to estimate $ \big| \big| L A (\alpha A+i)^{-N}  u \big| \big|_{L^2} $, when $ L = D_r $ or $ r^{-1}|D_{\mathcal S}| $. We observe that, for $ k = 0 , 1 $,
$$ L A^k (\alpha A+i)^{-N} = (A^k- ki)\big( \alpha A + i (1-\alpha) \big)^{-N} L  , $$
which follows easily from (\ref{A6}), (\ref{A7}) and (\ref{formulaA}) (see also  (\ref{traverseresolventeA}) below). Thus, using (\ref{interpolationA}),
\begin{eqnarray*}
 \big| \big| L A (\alpha A+i)^{-N}  u \big| \big|_{L^2} & \leq & C || \big(\alpha A+i(1-\alpha) \big)^{-N} L u ||_{L^2}^{1-\frac{1}{N}} || L u ||_{L^2}^{\frac{1}{N}} \\
 & \leq & C \big| \big| L \big(\alpha A+i \big)^{-N}  u  \big| \big|_{L^2}^{1-\frac{1}{N}} \big| \big| L u  \big| \big|_{L^2}^{\frac{1}{N}} \\
 & \leq & C || \big(\alpha A+i \big)^{-N}  u ||_{H_0^1}^{1-\frac{1}{N}} || u ||_{H_0^1}^{\frac{1}{N}} .
\end{eqnarray*}
The result follows. \finpreuve




\bigskip

We next record the basic Hardy inequality. Recall that we assume $ n \geq 3 $.

\begin{prop}  For all $ u \in C_0^{\infty}  $, we have
\begin{eqnarray}
 || r^{-1} u ||_{L^2} \leq \frac{2}{n-2} || \partial_r u ||_{L^2}  . \label{Hardy1}
\end{eqnarray}
As a consequence, the multiplication by $ r^{-s} $, $ s \in [0,1] $, is bounded from $ H_0^{1} $ to $ L^2 $. Furthermore,
\begin{eqnarray}
 (\alpha A - \zeta)^{-1} r^{-s} u = r^{-s} (\alpha A- \zeta + i \alpha s)^{-1} u , \qquad (\alpha A - \zeta)^{-1} \partial_r u = \partial_r (\alpha A-\zeta+i\alpha)^{-1} u, \label{traverseresolventeA}
\end{eqnarray}
for all $ s \in [0,1] $, $ u \in H_0^1 $, $ \alpha > 0 $ and $ \zeta \in \Ca \setminus \Ra $ such that $ |\emph{Im}(\zeta)| > \alpha $.
\end{prop}

\noindent {\it Proof.} The inequality (\ref{Hardy1}) is a direct consequence of the same one dimensional inequality (on $ L^2 (\Ra^+ , r^{n-1}dr) $) which is standard. The boundedness of $ r^{-1}$ on $ H_0^1 $, hence of $ r^{-s} $ by interpolation,   is then straighforward. The identities in (\ref{traverseresolventeA}) follow easily from (\ref{A6}), (\ref{A7}) and (\ref{formulaA}). \finpreuve

\bigskip

We will need later the following proposition.

\begin{prop} \label{injectiondomaine} The space $ \scal{r}^{-1} H_0^1 $ is contained in $ \emph{Dom}(A) $. Furthermore, for any symbol $ \sigma $ of order $ -1 $ ({\it i.e.} $ |\sigma^{(k)} (r) |\lesssim \scal{r}^{-1-k} $), there exists $ C > 0 $ such that
$$ || A ( \sigma (r) v ) ||_{L^2} \leq C ||  v  ||_{H_0^1} , $$
for all $ v \in  H_0^1 $.
\end{prop}

\noindent {\it Proof.} Let $ u = \sigma (r) v $ with $ v \in H_0^1 $. Let $ v_k \in C_0^{\infty} $ be a sequence approaching $v$ in $ H_0^1 $. Then, for all $ w \in \mbox{Dom}(A) $, we have
$$ (A w , u )_{L^2} = \lim_k (A w , \sigma (r) v_k)_{L^2} = \lim_{k \rightarrow \infty} (w, A ( \sigma (r) v_k))_{L^2} = (w , B v )_{L^2} , $$ 
where $ B $ is the (closure to $ H_0^1 $ of the) differential operator
$$ r \sigma (r) D_r + \frac{n}{2i} \sigma(r) + \frac{1}{i} r \sigma^{\prime} (r) , $$
which is bounded on $ H_0^1 $.
In particular
$$ || A u ||_{L^2} = || B v ||_{L^2} \leq C \big( || D_r v ||_{L^2} + || v ||_{L^2} \big) \leq C ||  v ||_{H_0^1} , $$
and this completes the proof. \finpreuve

\bigskip

We finally record simple weighted estimates. When $ W $ is a function, we set
$$ || W \nabla_{g_0} u ||_{L^2}^2  = || W \partial_r u ||_{L^2}^2 + \big| \big| W r^{-1}|D_{\mathcal S}| u  \big| \big|_{L^2}^2 . $$

\begin{prop} \label{propositionalgebreanalyse} There exists $ C > 0 $ such that for all $ u , v \in C_0^{\infty}  $, all non vanishing smooth function $ W : (0,+\infty) \rightarrow \Ca $ of $r$, and all admissible perturbation $ K $
\begin{eqnarray}
\big| \big( u , \emph{div}_{g_0} (K^{\rm sc} dv) \big) \big| \leq C ||| K |||_0   ||W (r) \nabla_{g_0} u ||_{L^2}   || W(r)^{-1} \nabla_{g_0} v ||_{L^2}  . \label{algebreanalyse}
\end{eqnarray}
\end{prop}

\noindent {\it Proof.} We start by writing
$$ \big( u , \mbox{div}_{g_0} (K^{\rm sc} dv) \big) = \int_{\Ra^+ \times {\mathcal S}} \left( \begin{matrix} W\partial_r \overline{u} \\  W d_{\mathcal S} \overline{u}/r \end{matrix} \right) \cdot  \left( \begin{matrix} K_{11} & K_{21}^{\dag} \\ K_{21}  & K_{22} \end{matrix} \right) \left( \begin{matrix} W^{-1}\partial_r v \\  W^{-1} d_{\mathcal S} v /r \end{matrix} \right) d \mu , $$
since the multiplication by $ W $ commutes with $ K $. Then, using (\ref{Linfini1}), (\ref{Linfini2}) and the fact that, if $ \xi \in T_{\omega}^* {\mathcal S} , V \in T_{\omega} {\mathcal S} $, $ |\xi \cdot V| \leq |\xi|_{h_{0,\omega}^*} |V|_{h_{0,\omega}} $, one sees that $ \big| \big( u , \mbox{div}_{g_0} (K^{\rm sc} dv) \big) \big| $ is not greater than 
$$ |||K|||_0 \int_{(0,+\infty) \times {\mathcal S}} \big( |W \partial_r \overline{u}|+  |W r^{-1} d_{\mathcal S} \overline{u}|_{h_{0,\omega}^*} \big) \big( |W^{-1} \partial_r v| + |W^{-1} r^{-1} d_{\mathcal S} v|_{h_{0,\omega}^*} \big) d \mu  . $$
By using the Cauchy-Schwarz inequality combined with the fact that
\begin{eqnarray*}
 \int_{\mathcal S} |W^{-1} r^{-1} d_{\mathcal S} v|_{h_{0,\omega}^*}^2 d {\rm vol}_{h_0} = r^{-2}W^{-2} \int_{\mathcal S} | d_{\mathcal S}  v|_{h_{0,\omega}^*}^2 d {\rm vol}_{h_0} & = & r^{-2} W^{-2} \big( v , -\Delta_{\mathcal S}  v \big)_{L^2 ({\mathcal S})}  \\
 & = & \big| \big|  W^{-1}(r) r^{-1} |D_{\mathcal S}| v (r) \big| \big|_{L^2 (\mathcal S)}^2
\end{eqnarray*} 
 the conclusion follows easily. \finpreuve
 
\bigskip

\subsection{Jensen-Mourre-Perry estimates} \label{soussectionMourre}
In this subsection, we define closed realizations of  operators of the form $ P_K $ (see (\ref{definitionPK})) and prove resolvent estimates thereon. Although we follow closely the Jensen-Mourre-Perry techniques \cite{JMP}, the proofs of resolvent estimates will be self contained. 
We have to review the proof for we will need to control most  estimates with respect to the perturbation $K$ and also since we need to prove $ H^{-1} \rightarrow H_0^1 $ estimates. 

Starting from (\ref{definitionformeK}), we observe first that the sesquilinear form $ Q_K $ satisfies
$$ |Q_K (u,v)| \leq C || u ||_{H_0^1} || v ||_{H_0^1} , $$
for all $ u , v \in C_0^{\infty} $ hence has a unique continuous extension to $ H_0^1 \times H_0^1 $. Everywhere in the sequel, we denote this extension by $ \bar{Q}_K $.

\begin{prop} \label{pourfirstitem} Let $ K $ be a $N$-admissible perturbation tensor.
\begin{itemize}
\item{The operator $ P_K : C_0^{\infty} \rightarrow C_0^{\infty} $ has a unique linear continuous extension $ \bar{P}_K :  H_0^1 \rightarrow H^{-1} $ and this extension satisfies
$$ ( \bar{P}_K u , v )  = \bar{Q}_K (u,v) , \qquad u , v \in H_0^1 . $$ Furthermore, there exists $ C $ independent of $K$ such that
\begin{eqnarray}
\big( 1 - C ||| K |||_0 \big) || u ||_{\dot{H}_0^1}^2 \leq (\bar{P}_K u , u ) \leq \big( 1 + C ||| K |||_0 \big) || u ||_{\dot{H}_0^1}^2 , \label{minorationhomogene}
\end{eqnarray}
for all $ u \in H_0^1 $.}
\item{If in addition $|| K ||_{L^{\infty}}$ is small enough, then one defines a self-adjoint operator $ \hat{P}_K $ on $ L^2 $ by
\begin{enumerate}
\item{$ \emph{Dom} (\hat{P}_K) = \{ u \in H_0^1 \ | \ \mbox{there exists} \ C_u > 0 \  |\bar{Q}_K (u,v)| \leq C_u || v ||_{L^2}, \ \ \mbox{for all} \ v \in H_0^1   \} $.}
\item{If $ u \in \emph{Dom} (\hat{P}_K) $, $ \hat{P}_K u  $ is defined as the unique element of $ L^2 $ such that 
$$ ( \hat{P}_K u , v )_{L^2} = \bar{Q}_K (u,v) , \qquad v \in H_0^1 . $$}
\item{The operator $ \hat{P}_K $ is nonnegative and has the property that
$$ \emph{Dom} (\hat{P}_K^{1/2}) = H_0^1, \qquad || u ||_{H_0^1}/2 \leq || (1+\hat{P}_K)^{1/2} u ||_{L^2} \leq 2 || u ||_{H_0^1} . $$}
\end{enumerate}}
\end{itemize}
\end{prop}


The proof is standard hence is omitted. We simply note that (\ref{minorationhomogene}) follows straightforwardly by density from (\ref{definitionPK}) and Proposition \ref{propositionalgebreanalyse} (with $ W \equiv 1 $), once noticed that $ (\bar{P}_0 u , u ) = || u ||_{\dot{H}_0^1}^2 $.

\bigskip





We now study the resolvent of $ \hat{P}_K $. Let us introduce the notation
$$ \Sigma (\varepsilon_0,\delta_0) = \left\{ (\delta, \varepsilon) \in \Ra \times \Ra \ | \ 0 < |\delta| \leq \delta_0, \ \ |\varepsilon| \leq \varepsilon_0 , \ \varepsilon \delta \geq 0 \right\} . $$
We also set
\begin{eqnarray}
 P_K (\varepsilon) = P_K + \sum_{k=1}^N \frac{\varepsilon^k}{k!} \mbox{ad}_A^k P_K  , \label{pourlepremierterme}
\end{eqnarray}
where, as usual, $ \mbox{ad}_A P = [P,A] = PA - A P $ and $  \mbox{ad}_A^{k+1} P = [\mbox{ad}_A^k  P, A] $. This definition makes sense as an equality between operators on $ C_0^{\infty} $. Its interest is the following easily verified property,
\begin{eqnarray}
 \big[ A , P_K (\varepsilon) \big] = - \frac{\partial}{\partial \varepsilon} P_K (\varepsilon) - \frac{\varepsilon^N}{N!} \mbox{ad}_A^{N+1} P_K . 
 \label{calculdederivee}
\end{eqnarray}
Using (\ref{commutateurformel0}), we have the formula
\begin{eqnarray}
 \mbox{ad}_A^k P_K = \frac{2^k}{i^k} P_{K_k}, \qquad K_k := \left(1 - \frac{r \partial_r}{2} \right)^k K . \label{Kkexplicite}
\end{eqnarray}
We can thus rewrite (\ref{pourlepremierterme}) as
\begin{eqnarray}
P_K (\varepsilon) = \left(P_K + \sum_{2 \leq 2 j \leq N} (-1)^j\frac{(2\varepsilon)^{2j}}{(2j)!} P_{K_{2j}} \right) - 2 i \varepsilon \left(  P_{K_1} + \sum_{3 \leq 2 j + 1 \leq N} (-1)^j \frac{(2\varepsilon)^{2j}}{(2j+1)!} P_{K_{2j+1}}  \right) \label{aussisymmetrie}
\end{eqnarray}
where both brackets are symmetric on $ C_0^{\infty} $ with respect to $ d \mu $. The operator $ P_K (\varepsilon) $ can be extended to an $ H_0^1 \rightarrow H^{-1} $ operator by Proposition \ref{pourfirstitem} (item 1), ie
\begin{eqnarray}
 \bar{P}_K (\varepsilon) = \bar{P}_K + \sum_{k=1}^N \frac{(-2 i\varepsilon)^k}{k!}  \bar{P}_{K_k} . \label{polynomeepsilon} 
\end{eqnarray} 
Similarly, the identity (\ref{aussisymmetrie}) can accordingly be extended as an equality between  $ H_0^1 \rightarrow H^{-1} $ operators.  




\begin{prop} \label{proplongue} Fix $ N \geq 1 $ and $ C > 0 $. There exist $ \varrho > 0 $ and $ \delta_0, \varepsilon_0 > 0 $ such that, for all $ N $-admissible perturbations satisfying
$$ ||| K |||_0 + ||| K |||_1 \leq \varrho, \qquad |||K|||_2 + \cdots + ||| K|||_N \leq C  $$
and all $ (\delta,\varepsilon) \in \Sigma (\varepsilon_0 , \delta_0) $,
we have the following results.
\begin{enumerate}
\item{The operator
$$ \bar{P}_K (\varepsilon) - (1+i\delta) \bar{I} : H_0^1 \rightarrow H^{-1} $$
is a bounded isomorphism (see Proposition \ref{injection} for $ \bar{I} $). Its inverse
$$ R (\varepsilon,\delta):= \big( \bar{P}_K (\varepsilon) - (1+i\delta) \bar{I} \big)^{-1}  $$ 
satisfies
\begin{eqnarray}
 \left| \left| R (\varepsilon , \delta) \right| \right|_{H^{-1} \rightarrow H_0^1} \leq \min\left( \frac{10}{|\varepsilon|} , \frac{8}{|\delta|} \right) . \label{borneepsilon} 
 \end{eqnarray}
}
\item{When $ \varepsilon = 0 $ and $ f \in L^2 $,
$$ \big( \bar{P}_K (0) - (1+i \delta) \bar{I} \big)^{-1} \bar{I} f = (\hat{P}_K - 1-i\delta)^{-1} f . $$}
\item{For all $ (\varepsilon,\delta) \in \Sigma (\varepsilon_0,\delta_0) $,
$$ R (\varepsilon , \delta)^* = R (- \varepsilon , - \delta) . $$}
\end{enumerate}
\end{prop}

\noindent {\it Proof.} By symmetry of $K$, hence of all $ K_k $,  $ ( \bar{P}_{K_k} u , u ) $ is real  for all $ u \in H_0^1 $. 
Using that $ ||| K |||_0 $ is small enough, it follows from  (\ref{minorationhomogene}) and (\ref{aussisymmetrie}) that
$$ \left( \frac{3}{4} - \gamma \varepsilon^2 \right) ||  u ||_{\dot{H}_0^1}^2 - || u ||_{L^2}^2  \leq \mbox{Re} \left(  \big(\bar{P}_{K}(\varepsilon) - (1+i \delta) \bar{I} \big) u , u \right) \leq \left( \frac{5}{4} + \gamma \varepsilon^2 \right) ||  u ||_{\dot{H}_0^1}^2 - || u ||_{L^2}^2 , $$
for some constant $ \gamma > 0 $ independent of $K $ as long as $ |||K|||_0 + \cdots + ||| K|||_N $ remains bounded.
 Similarly, if $ |||K_1|||_0 $ is small enough,  we have
$$ \mbox{sgn}(\varepsilon) \mbox{Im} \left(  \big(\bar{P}_{K}(\varepsilon) - (1+i \delta)\bar{I} \big) u , u \right) \geq |\varepsilon|  (1 - \gamma \varepsilon^2) ||  u ||_{\dot{H}_0^1}^2 + |\delta|  \big| \big| u  \big| \big|_{L^2}^2 $$
since the first term of $ \bar{P}_K (\varepsilon) $ contributing in the imaginary part is $ 2 \varepsilon ( \bar{P}_{K_1} u , u) $ and  $ \mbox{sgn}(\varepsilon) \delta = |\delta| $ by definition of $ \Sigma(\varepsilon_0, \delta_0) $. 
If $ |\varepsilon| \leq \varepsilon_0 $ is small enough, we   obtain
\begin{eqnarray}
 \mbox{Re} \left(  \big(\bar{P}_{K}(\varepsilon) - (1+i \delta) \bar{I} \big) u , u \right) & \leq & \frac{3}{2}  ||  u ||_{\dot{H}_0^1}^2 - || u ||_{L^2}^2 ,
 \label{minoration1ajoute} \\
 \mbox{Re} \left(  \big(\bar{P}_{K}(\varepsilon) - (1+i \delta) \bar{I} \big) u , u \right) & \geq & \frac{1}{2}  ||  u ||_{\dot{H}_0^1}^2 - || u ||_{L^2}^2 ,
 \label{minoration1} \\
\mbox{sgn}(\varepsilon) \mbox{Im} \left(  \big(\bar{P}_{K}(\varepsilon) - (1+i \delta) \bar{I} \big) u , u \right) & \geq & \frac{|\varepsilon|}{2} || u ||_{\dot{H}_0^1}^2 + |\delta| || u ||_{L^2}^2 . \label{minoration2}
\end{eqnarray}
We wish to get lower bounds in term of the $ H_0^1 $ norm.
Let $ \theta = \theta (\varepsilon,\delta) $  be such that
$ \cos (\theta) = - |\varepsilon| /4 $ and $ \mbox{sgn}(\varepsilon) \sin (\theta) = (1- \varepsilon^2/16)^{1/2} $.
 We  then have the coercivity estimate
\begin{eqnarray}
 \mbox{Re} \left( e^{i \theta}  \big(\bar{P}_{K}(\varepsilon) - (1+i \delta) \bar{I} \big) u , u \right)  \geq \frac{|\varepsilon|}{10} || u ||_{H_0^1}^2 , 
 \label{coercivite}
\end{eqnarray}
where we now have the $ H_0^1 $ norm in the right-hand side. Indeed, using that the duality $ (.,.) $ is antilinear in the first factor, the left hand side of (\ref{coercivite}) reads
$$  \cos \theta  \mbox{Re} \left( \big(\bar{P}_{K}(\varepsilon) - (1+i \delta) \bar{I} \big) u , u \right) + \sin \theta  \mbox{Im} \left( \big(\bar{P}_{K}(\varepsilon) - (1+i \delta) \bar{I} \big) u , u \right) . $$
Multiplying (\ref{minoration1ajoute}) by $ \cos \theta $ (which is negative) and (\ref{minoration2}) by $ \mbox{sgn}(\varepsilon) \sin (\theta) $ allows to bound from below this expression  by
$$ \frac{|\varepsilon|}{4} || u ||_{L^2}^2 + |\varepsilon| \left( \frac{(1-\varepsilon^2/16)^{1/2}}{2} - \frac{3}{8}\right) || u ||_{\dot{H}_0^1}^2 \geq \frac{|\varepsilon|}{10} ||u||_{H_0^1}^2 $$
 if $ \varepsilon $ is small enough. It follows from (\ref{coercivite}) that the operator $ \bar{P}_{K}(\varepsilon) - (1+i \delta) \bar{I} $ is injective and has a closed range (this would also follow from (\ref{minoration2})). Using the usual Lax-Milgram argument, the estimate (\ref{coercivite}) implies that $ e^{i \theta} \big(\bar{P}_{K}(\varepsilon) - (1+i \delta) \bar{I} \big) $ and hence $\bar{P}_{K}(\varepsilon) - (1+i \delta) \bar{I} $  are isomorphisms between $ H_0^1 $ and $ H^{-1} $,
%
%
which proves the existence of $ R (\varepsilon,\delta) $. To complete the proof of the first item, it remains to prove (\ref{borneepsilon}). The bound  $ 10 / |\varepsilon| $ follows from (\ref{coercivite}). We prove the bound  $  8/ |\delta|$ in a similar fashion as follows. We choose $ \beta \in \Ra $ such that  $ \mbox{sgn}(\varepsilon) \sin (\beta) = (1-\delta^2/4)^{1/2} $  and 
$ \cos (\beta) = |\delta|/2 $. Then, it is not hard to see as above that using (\ref{minoration1}) and (\ref{minoration2}) we have
$$ \mbox{Re} \left( e^{i \beta} \big(\bar{P}_{K}(\varepsilon) - (1+i \delta) \bar{I} \big) u , u \right)  \geq \frac{|\delta|}{8} || u ||_{H_0^1}^2 , $$
provided that $ \delta $ is small enough. This last estimate and (\ref{coercivite}) imply (\ref{borneepsilon}).

To prove the second item, we observe that $ \big( \bar{P}_K (0) - (1+i \delta) \bar{I} \big)^{-1} \bar{I} f $ is the unique $ u \in H_0^1 $ such that
$$ \bar{Q}_K (u,v) - \big( (1+i \delta) u , v \big)_{L^2} = (f,v)_{L^2} , $$
for all $v \in H_0^1 $. This implies precisely that $ u $ belongs to $ \mbox{Dom} (\hat{P}_K) $ and that 
$$ \big( \hat{P}_K u - (1+ i \delta) u , v \big)_{L^2} = (f,v)_{L^2} , $$
which shows that $ u = (\hat{P}_K - 1 -i \delta)^{-1} f $. 

 To prove the third item, we use the definition (\ref{symetriser}) to
write
\begin{eqnarray}
 ( f,R(-\varepsilon,-\delta) g ) & = & \left( \big(\bar{P}_{K}(\varepsilon) - (1+i \delta) \bar{I} \big) R (\varepsilon,\delta) f , R (-\varepsilon,-\delta) g \right) \nonumber \\
 & = &  \left( R (\varepsilon,\delta) f , \big(\bar{P}_{K}(-\varepsilon) - (1-i \delta) \bar{I} \big)  R (-\varepsilon,-\delta) g \right) \nonumber \\
 & = & ( R(\varepsilon,\delta)f,g ) , \nonumber
\end{eqnarray} 
since, by the definition of $ \bar{P}_K (\pm \varepsilon) $, it is clear that $ (\bar{P}_K(\varepsilon) u , v ) = ( u,\bar{P}_K(-\varepsilon)  v ) $ for all $ u , v \in H_0^1 $ (see for instance (\ref{aussisymmetrie})). The result follows.
\finpreuve

\bigskip
We next recall a classical lemma (see \cite{JMP}) on differential inequalities.

\begin{lemm} \label{inegalitedifferentielle} Let $ C > 0 $, $ \varepsilon_0 > 0 $, $ \gamma > 0 $ and $ 0 \leq  \beta < 1 $  be  fixed constants. Then there exists $ C^{\prime} > 0 $ such that, for all  differentiable map $ F : (0,\varepsilon_0) \rightarrow {\mathcal L}(H^{-1},H_0^1) $   satisfying
\begin{eqnarray}
 \left| \left| \frac{d}{d \varepsilon} F (\varepsilon) \right| \right|_{H^{-1} \rightarrow H_0^1} & \leq & C ( || F (\varepsilon) ||_{H^{-1} \rightarrow H_0^{1}} + 1) \varepsilon^{- \beta} ,  \label{surladerivee} \\
  || F (\varepsilon) ||_{H^{-1} \rightarrow H_0^1} & \leq & C \varepsilon^{-\gamma} , \label{surlafonction}
\end{eqnarray}
for all $ \varepsilon \in (0,\varepsilon_0) $, we have
$$ || F (\varepsilon) ||_{H^{-1} \rightarrow H_0^1} \leq C^{\prime}, $$ for all $ \varepsilon \in (0,\varepsilon_0) $.
\end{lemm}

\noindent {\it Proof.} Consider the sequence $ (\gamma_k)_{k \in \Na} $ defined by
$$ \gamma_0= \gamma, \qquad \gamma_{k+1} = \begin{cases} \gamma_k + \beta - 1 & \mbox{if} \ \gamma_k + \beta > 1 , \\
0 & \mbox{if} \ \gamma_k + \beta \leq 1 . \end{cases} $$
It is easy to check that $ \gamma_k = 0 $ for all $k$ large enough. The lemma then follows from the observation that for all $ k \geq 0 $ there exists $ C_k > 0 $ such that, for all $ F $ satisfying (\ref{surladerivee}) and (\ref{surlafonction}), 
$$ || F (\varepsilon) ||_{H^{-1} \rightarrow H_0^1} \leq C_k \varepsilon^{-\gamma_k} , \qquad \varepsilon \in (0,\varepsilon_0) , $$
which is obtained by an elementary induction. \finpreuve


\bigskip


In the following proposition $ K $ is a fixed $ N+1 $-admissible perturbation. For simplicity, when $ k \geq 1 $ is an integer, we will use $ R (\varepsilon,\delta)^k $ as the obvious short hand 
for $ ( R (\varepsilon,\delta) \bar{I} )^{k-1} R (\varepsilon,\delta) $.
\begin{prop} \label{differentiationfaible} For all $ (\varepsilon,\delta) \in \Sigma(\varepsilon_0,\delta_0) $, the function $ t \mapsto e^{itA} R (\varepsilon,\delta) e^{-itA} $ can be weakly differentiated at $ t = 0 $ and 
\begin{eqnarray}
 \frac{1}{i} \left. \frac{d}{dt} \big( e^{itA} R (\varepsilon,\delta) e^{-itA} \big) \right|_{t=0} = - \frac{\partial}{ \partial \varepsilon} R (\varepsilon,\delta) + R (\varepsilon,\delta) \left( (-2i)^{N+1} \frac{\varepsilon^N}{N!} \bar{P}_{K_{N+1}} \right) R (\varepsilon,\delta) . 
 \label{thefirstassertion}
\end{eqnarray}
Similarly, $\frac{1}{i} \left. \frac{d}{dt} \big( e^{itA} R (\varepsilon,\delta)^N e^{-itA} \big) \right|_{t=0} $, defined in the weak sense, reads
$$   - \sum_{k=0}^{N-1} R(\varepsilon,\delta)^k \left( \frac{\partial}{ \partial \varepsilon} R (\varepsilon,\delta) - R (\varepsilon,\delta) \left( (-2i)^{N+1} \frac{\varepsilon^N}{N!} \bar{P}_{K_{N+1}} \right) R (\varepsilon,\delta) \right) R (\varepsilon,\delta)^{N-1-k} . $$
\end{prop}

Recall that the action of $ e^{itA} $ on $ H_0^1 $ and $ H^{-1} $ respectively is described in Proposition \ref{differentiabiliteH01}. We also point out that the derivative $ \partial R (\varepsilon,\delta) / \partial \varepsilon $ is well defined, in the $ {\mathcal L}(H^{-1},H_0^1) $ topology, since $ R (\varepsilon,\delta) $ is the inverse of $ \bar{P}_K (\varepsilon) - (1 + i \delta) \bar{I} $ which depends polynomially on $ \varepsilon $  in $ {\mathcal L}(H_0^1,H^{-1}) $ (see (\ref{polynomeepsilon})).

\bigskip

\noindent {\it Proof.} Since $ e^{itA} $ is an isomorphism on $ H_0^1 $ and $ e^{-itA} $  an isomorphism on $ H^{-1} $, $ e^{itA} R (\varepsilon,\delta) e^{-itA} $ is the inverse of 
$$ e^{itA} \left( \bar{P}_K (\varepsilon)  - (1+i \delta) \bar{I} \right) e^{-itA} = e^{-2t} \bar{P}_{K^t} (\varepsilon) - (1+i\delta) \bar{I} , $$
where the equality with the right-hand side follows from (\ref{commutateurparderivee}) and the first item of Proposition \ref{pourfirstitem}.
We claim that this operator can be weakly differentiated in $t$. Indeed, if $ u , v \in H_0^1 $, we have
\begin{eqnarray}
 \big( e^{-2t} \bar{P}_{K^t}(\varepsilon) u , v \big) & = &  \big(\bar{P}_{K}(\varepsilon) u , v \big) -2 \int_0^t \big(\bar{P}_{K^s_1}(\varepsilon) u , v \big) ds . \label{moitie1}
\end{eqnarray}
This is easily seen first with $ u , v \in C_0^{\infty} $ by using (\ref{commutateurparderivee}), and then on $ H_0^1 $ by density.
By writing $(\bar{P}_{K^s_1}(\varepsilon) u , v ) = \bar{Q}_{K_1^s(\varepsilon)}(u,v) $, we see that this quantity depends continously on $s$, by the strong continuity on $ L^2 $ of $ s \mapsto  (r \partial_r )^k K  (e^s r ) $, for $ k \leq N+1 $. This implies on one hand that 
$ ( e^{-2t}\widetilde{P}_{K^t}(\varepsilon) u , v ) $ can be differentiated at $t=0$ and on the other hand that
$$ \big| \big| e^{-2t}\bar{P}_{K^t}(\varepsilon) - \bar{P}_{K}(\varepsilon) \big| \big|_{H_0^1 \rightarrow H^{-1}} \leq C |t| $$
hence that
\begin{eqnarray}
 \big| \big| \big( e^{-2t}\bar{P}_{K^t}(\varepsilon) - (1+i \delta) \bar{I} \big)^{-1} - R (\varepsilon,\delta) \big| \big|_{H^{-1} \rightarrow H_0^{1}} \leq C_{\varepsilon,\delta} |t| . \label{sousderivee}
\end{eqnarray}
Using the resolvent identity, this shows that
$$ \big( e^{-2t}\bar{P}_{K^t}(\varepsilon) - (1+i \delta) \bar{I} \big)^{-1} - R (\varepsilon,\delta) = - R (\varepsilon,\delta) \big( e^{-2t}\bar{P}_{K^t}(\varepsilon) -\bar{P}_{K}(\varepsilon) \big) R (\varepsilon,\delta) + O (t^2) ,$$
where the $ O (t^2) $ holds in operator norm. This justifies the weak differentiability and the fact that
$$ \left. \frac{d}{dt} \big( e^{itA} R (\varepsilon,\delta) e^{-itA} \big) \right|_{t=0} = -  R (\varepsilon,\delta) \frac{d}{dt} \big( e^{itA} \bar{P}_K (\varepsilon) e^{-itA}\big)_{|t=0} R (\varepsilon,\delta) . $$
To compute the derivative, we use on one hand  (\ref{moitie1}) to see that
$$ \frac{d}{dt} \big( e^{itA}\bar{P}_{K}(\varepsilon) e^{-itA} \big)_{|t=0}  = - 2 \bar{P}_{K_1} (\varepsilon) . $$
On the other hand, using (\ref{polynomeepsilon}), we can compute $ \partial \bar{P}_K (\varepsilon) / \partial \varepsilon $ directly and check that
$$ - 2 \bar{P}_{K_1} (\varepsilon)  = - i\frac{\partial}{\partial \varepsilon} \bar{P}_K (\varepsilon)  -i  (-2i)^{N+1} \frac{\varepsilon^N}{N!} 
\bar{P}_{K_{N+1}} , $$
from which (\ref{thefirstassertion}) follows. In a similar fashion, we have
$$ e^{itA} R(\varepsilon,\delta)^N e^{-itA} =  R(\varepsilon,\delta)^N  + \sum_{k=0}^{N-1} R (\varepsilon,\delta)^k \left( e^{itA} R(\varepsilon,\delta)^N e^{-itA} - R (\varepsilon,\delta) \right) R (\varepsilon,\delta)^{N-1-k}  +O (t^2). $$
This allows to justify the weak differentiability and obtain the second assertion. \finpreuve

\bigskip

The next lemma is a convenient version of the standard quadratic estimates of \cite{Mourre}.

\begin{lemm} \label{lemmequadratique} Let $ B : H^{-1} \rightarrow H^{-1} $ be a bounded linear operator. Then for all $ (\varepsilon,\delta) \in \Sigma(\varepsilon_0,\delta_0) $ and all $K$ as in Proposition \ref{proplongue},
$$ \big| \big| R (\varepsilon,\delta) B \big| \big|_{H^{-1} \rightarrow H_0^1} \leq \frac{4}{|\varepsilon|^{1/2}} \big| \big| B^* R (\varepsilon,\delta) B \big| \big|_{H^{-1} \rightarrow H_0^{1}}^{1/2} . $$
\end{lemm}

\noindent {\it Proof.} Using (\ref{coercivite}), we have
$$ \mbox{Re} \left( e^{i \theta} \big(\bar{P}_{K}(\varepsilon) - (1+i \delta) \bar{I} \big) R (\varepsilon,\delta) B f , R (\varepsilon,\delta) B f \right)  \geq \frac{|\varepsilon|}{10} || R (\varepsilon,\delta) B f ||_{H_0^1}^2 , 
  $$
for all $ f \in H^{-1} $. The left-hand side reads
$$  \mbox{Re} \left( e^{i \theta}  B f , R (\varepsilon,\delta) B f \right) = \mbox{Re} \left( e^{i \theta}  f , B^* R (\varepsilon,\delta) B f \right) . $$
By bounding $ 1/10 $ from below by $ 1/16 $, this implies that
$$ || f ||_{H^{-1}}^2 || B^* R (\varepsilon,\delta) B  ||_{H^{-1} \rightarrow H^1_0} \geq \frac{|\varepsilon|}{16} || R (\varepsilon,\delta) B f ||_{H_0^1}^2 , $$
from which the result follows. \finpreuve

\bigskip

\begin{prop} \label{Mourrenonblackbox} Fix $ N \geq 1 $, $ M \geq 1 $ and $ 0 < \alpha < 1 $. There exist $ C > 0 $ large enough and $ \varrho, \varepsilon_0 , \delta_0 > 0 $ small enough such that, for all $ (\varepsilon,\delta) \in \Sigma (\varepsilon_0, \delta_0) $ and all $ N +1 $-admissible perturbation $ K $ satisfying
$$ ||| K |||_0 + ||| K |||_1 < \varrho, \qquad ||| K |||_2 + \cdots + ||| K |||_{N+1} \leq M , $$
we have
$$ \big| \big| (\alpha A+i)^{-N} R (\varepsilon, \delta)^N (\alpha A-i)^{-N} \big| \big|_{H^{-1} \rightarrow H_0^1} \leq C . $$
\end{prop}



\noindent {\it Proof.} Let us set first $ F^1 (\varepsilon) := (\alpha A-i)^{-1} R (\varepsilon,\delta) (\alpha A+i)^{-1} $. By  Lemma \ref{lemmequadratique} and the item 3 of Proposition \ref{differentiabiliteH01} with $ \zeta = - i $, we have
\begin{eqnarray}
 || R (\varepsilon,\delta) (\alpha A+i)^{-1} ||_{H^{-1} \rightarrow H_0^1} \leq \frac{4}{|\varepsilon|^{1/2}} || F^1 (\varepsilon)  ||^{1/2}_{H^{-1} \rightarrow H_0^1} . \label{bornebientot}
\end{eqnarray}
By taking the adjoint and using the  third items of Propositions \ref{differentiabiliteH01} and \ref{proplongue},  the same estimate holds for $  || (\alpha A-i)^{-1} R (\varepsilon,\delta)  ||_{H^{-1} \rightarrow H_0^1} $.
On the other hand, by using the item 4 of Proposition \ref{differentiabiliteH01} and Proposition \ref{differentiationfaible}, we obtain
\begin{eqnarray*}
 \frac{d}{d \varepsilon} F^1 (\varepsilon) & = & i A (\alpha A-i)^{-1} R (\varepsilon,\delta) (\alpha A+i)^{-1} - i (\alpha A-i)^{-1} R (\varepsilon,\delta) A(\alpha A+i)^{-1} + \\
 & &  (\alpha A - i)^{-1} R (\varepsilon,\delta) \left( (-2i)^{N+1} \frac{\varepsilon^N}{N!} \bar{P}_{K_{N+1}} \right) R (\varepsilon,\delta) (\alpha A+i)^{-1} , 
\end{eqnarray*}
as an equality between $ H^{-1} \rightarrow H_0^1 $ operators.
Therefore, using (\ref{bornebientot}) and the bound (\ref{borneepsilon}) to handle the last term above, we get
$$ \left| \left| \frac{d}{d \varepsilon} F^1 (\varepsilon) \right| \right|_{H^{-1} \rightarrow H_0^1} \leq C_{\alpha} |\varepsilon|^{-1/2} || F^1 (\varepsilon) ||^{1/2}_{H^{-1} \rightarrow H_0^1} . $$
By (\ref{borneepsilon}) and Lemma \ref{inegalitedifferentielle}, we obtain that $ || F^1 (\varepsilon) ||_{H^{-1} \rightarrow H_0^1} \leq C $. In particular, the right-hand side of (\ref{bornebientot})
is at most of order $ |\varepsilon|^{-1/2} $. If we now set $ F^N (\varepsilon) := (\alpha A-i)^{-N} R (\varepsilon,\delta)^N (\alpha A+i)^{-N} $, we obtain similarly (using now the second part of Proposition \ref{differentiationfaible}) that
\begin{eqnarray}
 \left| \left| \frac{d}{d \varepsilon} F^N (\varepsilon) \right| \right|_{H^{-1} \rightarrow H_0^1} & \lesssim_{\alpha } & \! \! \! \! \! \big|\big| A (\alpha A +i)^{-N} R (\varepsilon,\delta)^N (\alpha A-i)^{-N} \big|\big|_{H^{-1} \rightarrow H_0^1} + \nonumber \\ & & \! \! \! \! \! \big|\big|  (\alpha A +i)^{-N} R (\varepsilon,\delta)^N (\alpha A-i)^{-N} A \big|\big|_{H^{-1} \rightarrow H_0^1} + \nonumber
 \\ &  & \! \! \! \! \! \! |\varepsilon| \big| \big| (\alpha A+i)^{-N} R (\varepsilon,\delta) \big| \big|_{H^{-1} \rightarrow H_0^1}  \big| \big|  R (\varepsilon,\delta) (\alpha A -i)^{-N} \big| \big|_{H^{-1} \rightarrow H_0^1} .
 \nonumber
\end{eqnarray} 
The last line is the contribution of
$$ \sum_{k=0}^{N-1} (\alpha A + i)^{-N} R (\varepsilon,\delta) \left\{  R (\varepsilon,\delta)^k \varepsilon^N \bar{P}_{K_{N+1}} R (\varepsilon,\delta)^{N-1-k}  \right\} R (\varepsilon,\delta) (\alpha A - i)^{-N} , $$
where each bracket $ \{ \cdots \} $ in the middle has a $ H^{-1} \rightarrow H^{-1} $ norm of order $ \varepsilon $ by (\ref{borneepsilon}). 
This is then bounded since the right-hand side of (\ref{bornebientot}) is at most of order $ |\varepsilon|^{-1/2} $. Then,  estimating $  ||R(\varepsilon,\delta)^N (\alpha A-i)^{-N} ||_{H^{-1} \rightarrow H_0^1} $ by
\begin{eqnarray*}
 ||R(\varepsilon,\delta) ||_{H^{-1} \rightarrow H_0^1}^{N-1}
 ||R(\varepsilon,\delta) (\alpha A-i)^{-1} ||_{H^{-1} \rightarrow H_0^1} || (\alpha A-i)^{1-N} ||_{H^{-1} \rightarrow H^{-1}}  \lesssim 
  |\varepsilon|^{-(N-1) - \frac{1}{2}} 
\end{eqnarray*}  which follows from (\ref{borneepsilon}) and (\ref{bornebientot}), and using a  similar estimate for $(\alpha A + i)^{-N} R (\varepsilon,\delta)^N $ together with the interpolation estimate of Proposition \ref{differentiabiliteH01}, we thus obtain
$$ \left| \left| \frac{d}{d \varepsilon} F^N (\varepsilon) \right| \right|_{H^{-1} \rightarrow H_0^1}  \lesssim_{\alpha } 1 + || F^N (\varepsilon) ||_{H^{-1} \rightarrow H_0^1}^{1 - \frac{1}{N}} \left( |\varepsilon|^{ \frac{1}{2}-N} \right)^{\frac{1}{N}}
$$
so the conclusion follows again from Lemma \ref{inegalitedifferentielle}. \finpreuve

\bigskip

We end up this section with two technical results which will be useful when we will ultimately replace the powers of $ (\alpha A \pm i)^{-1} $ by powers of $\scal{r}^{-1} $. More precisely, to prove Proposition \ref{parametrixe}, we will combine the estimates of Proposition \ref{Mourrenonblackbox} with refinements of Hardy type inequalities, as they appear for instance in \cite{VasyWunschHardy}, and which we now consider. 

\begin{lemm} \label{pour2eme} Let $ 0 \leq s < \frac{n-2}{2} $. Then, there exists $ C_s > 1 $ such that for all $ u \in C_0^{\infty} $ and all $ \delta \geq 0 $,
$$ 
C_s^{-1} || \partial_r \big( (r+\delta)^{-s}u \big) ||_{L^2  } \leq || (r+\delta)^{-s} \partial_r u ||_{L^2} \leq C_s || \partial_r \big((r+\delta)^{-s} u \big) ||_{L^2 } $$
and
$$ || (r+\delta)^{-s-1} u ||_{L^2} \leq C_s || (r+\delta)^{-s} \partial_r u ||_{L^2  }  .  $$
\end{lemm}

This last estimate is a Hardy inequality which is very close to \cite[Lemma 3.2]{VasyWunschHardy}. Here the additional information is the equivalence of the norms of $ (r+\delta)^{-s} \partial_r u $ and $ \partial_r ((r+\delta)^{-s} u) $ which we will need below.

\bigskip

\noindent {\it Proof.} By decomposing $u$ along an orthonormal basis of $ L^2 ({\mathcal S}, d{\rm vol}_{h_0}) $, it suffices to prove the result for functions $ v \in C_0^{\infty}(\Ra^+) \subset L^2 (\Ra^+ , r^{n-1}dr) $. Using the straighforward computation 
$$ \partial_r ((r+\delta)^{-s}v) - (r+\delta)^{-s} \partial_r v =  - s (r+\delta)^{-s-1} v, $$ and then using (\ref{Hardy1}), we obtain
$$ \left| || \partial_r( (r+\delta)^{-s}v) ||_{L^2 (r^{n-1}dr)} - || (r+\delta)^{-s} \partial_r v ||_{L^2 (r^{n-1}dr)} \right| \leq \frac{2s}{n-2} || \partial_r( (r+\delta)^{-s}v) ||_{L^2 (r^{n-1}dr)} . $$
Since $ 2s / (n-2) < 1 $, we obtain the equivalence of the norms. Using (\ref{Hardy1}), these norms control $ || (r+\delta)^{-s-1} v ||_{L^2(r^{n-1}dr)} $ and the result follows. \finpreuve

\bigskip

The next proposition is a generalization of an estimate which can be found in \cite[Proposition 4.1]{VasyWunschHardy}. We have to modify it to allow additional weights depending on $A$. In passing, we also get the full range of exponents $ s \in ( 0 , (n-2)/2 ) $.

\begin{prop} \label{Ptrav} Fix $ 0 < s < (n-2)/2 $ and constants $ M , N \in \Na $. Then there exist $ C > 0 $ large enough and $ \varrho > 0 $ small enough such that, for all $ \alpha $ small enough, all $ \delta \in (0,1) $, all $N$-admissible perturbations $ K $ such that 
$$ ||| K |||_0 \leq \varrho \qquad \mbox{and} \qquad  ||| K |||_{N} \leq M , $$ we have
\begin{eqnarray}
  \big| \big| (r+\delta)^{-s} \nabla_{g_0} u \big| \big|_{L^2} \leq C \big| \big| (\alpha A+i)^{-N} \hat{P}_K (\alpha A+i)^{N} u \big| \big|^{1/2}_{L^2} \big| \big| (r+\delta)^{-2s} u \big| \big|^{1/2}_{L^2} , \label{reresult}
\end{eqnarray} 
for all $u \in (\alpha A+i)^{-N} \emph{Dom}(\hat{P}_K)$.
\end{prop}

Before giving the detailed proof, we recall first the nice basic idea on which it rests when $N = 0 $ and $ \hat{P}_K = - \Delta_{g_0} $, {\it i.e.} when $ K = 0 $. 
We  compute $  \mbox{Re} \ (-\Delta_{g_0} u,(r+\delta)^{-2s} u)_{L^2} $ or more precisely $   \mbox{Re} \ \bar{Q}_0 ( u,(r+\delta)^{-2s} u) $.
The condition $ \delta > 0 $ guarantees that $ (r+\delta)^{-2s} $ is bounded on $ H_0^1 $ by the item 1 of Proposition \ref{differentiabiliteH01}. By density of $C_0^{\infty} $ in $ H_0^1 $, on which $  \bar{Q}_0 $ is continuous, we may assume that $ u \in C_0^{\infty} $. Then, by integrating by part one finds
\begin{eqnarray}
   \mbox{Re}\ \bar{Q}_0( u,(r+\delta)^{-2s} u) & = &   s \int_{{\mathcal M}_0} r^{n-2} (r+\delta)^{-2s-2}  \big((n-2-2s)r + (n-1) \delta \big) |u|^2 d r d {\rm vol}_{h_0} \nonumber \\
 & & \ \ +  \ || (r+\delta)^{-s} D_r u ||_{L^2}^2 + \big| \big| (r+\delta)^{-s} r^{-1}| D_{\mathcal S}| u \big| \big|_{L^2}^2 \nonumber \\
 & \geq & || (r+\delta)^{-s} \nabla_{g_0} u ||_{L^2}^2   , \label{borneexpliciteHardyh}
\end{eqnarray}
since $ (n-2-2s)r + (n-1) \delta \geq 0 $ (here we may go up to $ s = (n-2)/2 $). If $ u $ belongs to $ \mbox{Dom}(\hat{P}_0) $, then $   \bar{Q}_0( u,(r+\delta)^{-2s} u) =  ( \hat{P}_0 u,(r+\delta)^{-2s} u)_{L^2} $ and, using the Cauchy-Schwarz inequality, this yields the result with $ C = 1 $ when  $ K = 0 $ and $ N = 0 $. The general case will follow from this model by a perturbation argument.

\bigskip

\noindent {\it Proof of Proposition \ref{Ptrav}.} We start by computing the commutator $ \big[ (\alpha A-i)^N , \hat{P}_K \big] $ in the form sense:  for all $ \psi, \varphi \in C_0^{\infty} ({\mathcal M}_0) $, we have
\begin{eqnarray}
 Q_K \big( (\alpha A+i)^N \psi , \varphi \big) = Q_K \big(\psi , (\alpha A-i)^N \varphi \big) + \sum_{k=1}^N C_N^k (2i \alpha)^k Q_{K_k} \big( \psi , (\alpha A-i)^{N-k} \varphi \big) , \label{local}
\end{eqnarray}
where $ K_k $ is as in (\ref{Kkexplicite}). This follows from (\ref{commutateurformel0}) and a simple induction on $N$.  Then, up to considering the closures of the quadratic forms, the identity (\ref{local}) remains true if one only assumes that $ \psi \in H_0^1  $ and $\varphi \in H_0^1 $  satisfy $ (\alpha A+i)^N \psi \in H_0^1 $ and $ (\alpha A-i)^{N} \varphi \in H_0^1 $ (this can be proved as in Lemma 3.1 of \cite{BoucletCPDE}\footnote{for convenience, we recall that the idea is to set $ \psi = (\alpha A+i)^{-N} v  $ and to use on one hand that $  (\alpha A+i)^{-N}  = \lim_{\tau \rightarrow \infty}  \frac{1}{i (N-1)!} \int_0^{\tau} (-it)^{N-1} e^{-t} e^{it \alpha A}   dt $, where the 
 integral preserves $ C_0^{\infty} $, and on the other hand that one can approximate $ v \in H_0^1 $ by some $ C_0^{\infty} $ function}). 
Consequently, by choosing $\psi = u $ and $ \varphi = (\alpha A-i)^{-N} (r + \delta)^{-2s} u $, we find that
$$ \mbox{Re} \ \big( (\alpha A+i)^{-N} \hat{P}_K (\alpha A+i)^N u , (r+\delta)^{-2s} u \big)_{L^2} = \mbox{Re} \ \bar{Q}_K \big( (\alpha A+i)^N u , (\alpha A-i)^{-N} (r + \delta)^{-2s} u \big) $$
can be written as 
\begin{eqnarray}
 \mbox{Re} \left( \bar{Q}_K \big(u , (r + \delta)^{-2s} u \big) + \sum_{k=1}^N C_N^k (2i \alpha)^k \bar{Q}_{K_k} \big( u , (\alpha A-i)^{-k} (r+\delta)^{-2s} u \big) \right) .
 \label{sumof} 
\end{eqnarray}
Our goal is to bound this expression from below similarly to (\ref{borneexpliciteHardyh}).
To study the contribution of the first term of (\ref{sumof}), we use Lemma \ref{propositionalgebreanalyse} and the expression of $ P_K - (-\Delta_{g_0}) $ given by (\ref{definitionPK}) to obtain
\begin{eqnarray*}
\big| \mbox{Re} \ \bar{Q}_{0}( u,(r+\delta)^{-2s} u)_{L^2} -  \mbox{Re} \ \bar{Q}_K ( u,(r+\delta)^{-2s} u) \big|  \lesssim \qquad \qquad \qquad \qquad \qquad \qquad \\
\qquad \qquad \qquad \qquad \qquad  ||| K |||_0 || (r+\delta)^{-s} \nabla_{g_0} u ||_{L^2}   \big( || (r+\delta)^{-s} \nabla_{g_0} u ||_{L^2} + || (r+ \delta)^{-s-1} u ||_{L^2} \big) .
\end{eqnarray*}
By using Lemma \ref{pour2eme} (which imposes $ s < (n-2)/2 $) and (\ref{borneexpliciteHardyh}), we get
\begin{eqnarray}
\mbox{Re} \ \bar{Q}_{K}( u,(r+\delta)^{-2s} u)_{L^2} \geq \big( 1 - C |||K|||_0 \big) || (r+\delta)^{-s} \nabla_{g_0} u ||_{L^2}^2 . \label{borneinf}
\end{eqnarray}
The result will then follow if we prove that, for each term in the sum of (\ref{sumof}),
\begin{eqnarray}
\big| Q_{K_k} \big( u , (\alpha A-i)^{-k} (r+\delta)^{-2s} u \big) \big| \leq C (1 + |||K_k|||_0) || (r+ \delta)^{-s} \nabla_{g_0} u ||^2_{L^2} . \label{ifweprove}
\end{eqnarray}
To justify (\ref{ifweprove}), we start by observing that for all $ v , w \in C_0^{\infty} $, we have 
\begin{eqnarray*}
 |Q_{K_k} (v,(r+\delta)^{-2s} w)| \leq C (1+|||K_k|||_0)  || (r+\delta)^{-s} \nabla_{g_0} v ||_{L^2}  || (r+\delta)^{s} \nabla_{g_0} (r+\delta)^{-2s} v ||_{L^2} .
\end{eqnarray*}
This follows easily from  of (\ref{algebreanalyse}) (strictly speaking (\ref{algebreanalyse}) only yields the contribution of $ K_k $ but the one of $ \Delta_{g_0} $ is similar). Using that
$$ || (r+\delta)^{s} \nabla_{g_0} (r+\delta)^{-2s} v ||_{L^2} \lesssim || (r+\delta)^{-s} \nabla_{g_0}  v ||_{L^2} + || (r+\delta)^{-s-1}  v ||_{L^2} , $$
and Lemma \ref{pour2eme}, and then a density argument to replace $ v $ and $w$ by any $ H_0^1 $ functions, we get
\begin{eqnarray*}
 |Q_{K_k} (u,(\alpha A-i)^{-k}(r+\delta)^{-2s} u)| \lesssim (1+|||K_k|||_0)  || (r+\delta)^{-s} \nabla_{g_0} u ||_{L^2}  || (r+\delta)^{-s} \nabla_{g_0}  \tilde{u} ||_{L^2} ,
\end{eqnarray*}
with
$$ \tilde{u} = (r+\delta)^{2s} (\alpha A-i)^{-k}(r+\delta)^{-2s} u . $$
The proof will then be complete if we show that
\begin{eqnarray}
 || (r+\delta)^{-s} \nabla_{g_0} \tilde{u}  ||_{L^2} \leq C || (r+\delta)^{-s} \nabla_{g_0} u  ||_{L^2} . \label{plusexplcitegradient}
\end{eqnarray} 
This is obtained as follows. Given $ \sigma \in \Ra $ and $ \zeta \in \Ca \setminus \Ra $, if $\alpha $ is small enough we have
\begin{eqnarray}
  || (r+\delta)^{\sigma} (\alpha A - \zeta)^{-k}(r+\delta)^{-\sigma}  ||_{L^2 \rightarrow L^2} \leq C  .\label{pourtraverserpoids}
\end{eqnarray}  
This comes from the first identity of (\ref{traverseresolventeA}), since $ (r+ \delta)^{|\sigma|} \approx r^{|\sigma|} + \delta^{|\sigma|} $ ({\it i.e.} their quotient is bounded from above and below). Therefore, using (\ref{traverseresolventeA}) and (\ref{pourtraverserpoids}) 
\begin{eqnarray*}
 || (r+\delta)^{-s} \partial_r  \tilde{u}||_{L^2} & \lesssim & || (r+\delta)^s \partial_r (\alpha A-i)^{-k} (r+\delta)^{-2s} u ||_{L^2} + || (r+\delta)^{s-1} (\alpha A-i)^{-k} (r+\delta)^{-2s} u ||_{L^2} \\
 & \lesssim & || (r+\delta)^s  (\alpha A-i + i \alpha)^{-k} \partial_r \big( (r+\delta)^{-2s} u \big) ||_{L^2} + ||  (r+\delta)^{-s-1} u ||_{L^2} \\
 & \lesssim & ||  (r+\delta)^{-s} \partial_r u  ||_{L^2} + ||  (r+\delta)^{-s-1} u ||_{L^2} \\
 & \lesssim & ||  (r+\delta)^{-s} \partial_r u  ||_{L^2}
\end{eqnarray*}
the last inequality following from Lemma \ref{pour2eme}. The estimate  
$$ \big| \big| (r+\delta)^{-s} r^{-1} |D_{\mathcal S}| \tilde{u} \big| \big|_{L^2} \lesssim \big| \big| (r+\delta)^{-s} r^{-1} |D_{\mathcal S}| u \big| \big|_{L^2} $$ is obtained in the very same way, using additionally that $ |D_{\mathcal S}| $ commutes with functions of $r$ and $A$.
This yields (\ref{plusexplcitegradient}) and completes the proof. \finpreuve

\section{Proofs of Propositions \ref{reductiongeometriqueprop} and \ref{parametrixe}} \label{Mourreonaconesuite}
\setcounter{equation}{0}
\subsection{Proof of Proposition \ref{reductiongeometriqueprop}}
Using (\ref{formedeG}), we recall that  $ g = a (r) dr^2 + 2 r b (r) dr + r^2 h (r) $ is a metric on $ (R_0,\infty) \times {\mathcal S} $ such that $ G = \kappa^* g $. This allows to recast the problem on a question on a half line times $ {\mathcal S} $: it suffices to show that one can find  $ \widetilde{R}_0 > 0 $ and a diffeomorphism
$$ \Xi : (\widetilde{R}_0 , \infty ) \times {\mathcal S} \rightarrow U \subset (R_0,\infty) \times {\mathcal S} $$
of the form
$$ \Xi (\tilde{r},\sigma) = \big( \bar{r} (\tilde{r},\sigma),\sigma \big) $$
such that
\begin{enumerate}
\item{for some symbol $ \xi \in  S^{-\rho} $ on $ (\widetilde{R}_0 , \infty ) \times {\mathcal S} $ ({\it i.e.} for all integer $k$
$ \partial_{\tilde{r}}^k \xi (\tilde{r},.) = {\mathcal O} (\tilde{r}^{-\rho-k}) $ in $ C^{\infty}({\mathcal S},\Ra) $)
$$ \bar{r} (\tilde{r},\sigma) =  \tilde{r} \big( 1 + \xi (\tilde{r},\sigma) \big) , $$}
\item{$ U $ contains $ (R_0^{\prime} , + \infty) \times {\mathcal S} $ for some $ R_0^{\prime} > R_0 $ large enough,  }
\item{at each point  $ (\tilde{r},\sigma) \in (\widetilde{R}_0 , \infty ) \times {\mathcal S} $, we have
\begin{eqnarray} 
 d {\rm vol}_{\Xi^* g} = \tilde{r}^{n-1} d \tilde{r} d {\rm vol}_{h_0} . \label{cequonveut}
\end{eqnarray} }
\end{enumerate}
To  build $ \bar{r} $ we check which conditions must be fulfilled. We first note that, at any point $ (r,\sigma ) \in (R_0,\infty) \times {\mathcal S} $,
$$ d {\rm vol}_{g} = F(r,\sigma) r^{n-1} d r d {\rm vol}_{h_0} , $$
with $ F - 1 \in S^{-\rho} $. Therefore, the condition (\ref{cequonveut}) reads
\begin{eqnarray}
 \tilde{r}^{n-1} = \frac{\partial \bar{r}}{\partial \tilde{r}}  \bar{r}^{n-1} F (\bar{r},\sigma) . \label{item3}
\end{eqnarray} 
If we assume that $ \Xi $ exists, the inverse diffeomorphism is of the form $ ( \underline{r} (r, \sigma) ,\sigma ) $. By evaluating (\ref{item3}) at $ \tilde{r} = \underline{r} (r,\sigma) $, we get
$$ \underline{r}^{n-1} \frac{1}{\frac{\partial \bar{r}}{\partial \tilde{r}}(\underline{r})} = r^{n-1} F (r,\sigma) , $$
that is
$$ \underline{r}(r,\sigma)^{n-1} \frac{\partial \underline{r}}{\partial r}(r,\sigma) = r^{n-1} F (r,\sigma) , $$
or equivalently
\begin{eqnarray}
 \frac{\partial \underline{r}^n}{\partial r} (r,\sigma) =  n r^{n-1} F (r, \sigma) , \label{equivalently} 
\end{eqnarray} 
which can be solved:
we write $ F $ as $ F (r,\sigma) = 1 + \delta (r,\sigma)  $ with $ \delta \in S^{-\rho} $ then, by following (\ref{equivalently}), we define for some $ R_1 > R_0 $
\begin{eqnarray}
  \underline{r} (r,\sigma) = \left( n \int_{R_1}^{r}  (1 + \delta(t,\sigma)) t^{n-1} dt \right)^{\frac{1}{n}} = r  \left( 1 - \frac{R_1^n}{r^n} + r^{-n} \int_{R_1}^{r}  \delta(t,q) t^{n-1} dt \right)^{\frac{1}{n}} , \label{pourinverser} 
\end{eqnarray} 
 for $ r  > R_1 $ and $ \sigma \in {\mathcal S} $. Since $ n \geq 2 $ and by assuming $ \rho \leq 1 $, it is not hard to check that
 $ r \mapsto \int_{R_1}^{r}  \delta(t,\sigma) t^{n-1} dt $
 belongs to $ S^{n-\rho} $ hence that the last bracket in (\ref{pourinverser}) is of the form $ 1 + S^{-\rho} $.
 It follows easily that, for $ R_2 \gg 1$ and for all $ \sigma \in {\mathcal S} $, $ r \mapsto \underline{r} (r,\sigma) $ is a diffeomorphism from $ (R_2 , \infty) $ to $ ( \underline{r} (R_2, \sigma),\infty ) $ hence that
 $$ (r,\sigma) \mapsto ( \underline{r} (r,\sigma),\sigma) $$
 is a diffeomorphism from $ (R_2,\infty) \times {\mathcal S} $ to an open subset containing $ \big( \sup_{\mathcal S}  \underline{r} (R_2,.) , \infty \big) \times {\mathcal S} $ which contains $ (\tilde{R}_0 , \infty ) \times {\mathcal S} $ for some $ \tilde{R}_0 \gg 1 $. The inverse diffeomorphism provides a diffeomorphism of the form $ \Xi :   (\tilde{r},\sigma) \mapsto ( \bar{r} (\tilde{r},\sigma) , \sigma ) $
 which, by construction, satisfies (\ref{item3}) and hence (\ref{cequonveut}).  Using that $ \underline{r} (r,\sigma) = r (1+S^{-\rho}) $ and by differentiating
 $ \underline{r} ( \bar{r} (\tilde{r},\sigma),\sigma) = \tilde{r} $, a routine analysis shows that
 $ \bar{r} (\tilde{r},\sigma) $ is of the form $ \tilde{r} (1 + S^{-\rho}) $,  which yields both item 1 and item 2. \finpreuve

\subsection{Proof of Proposition \ref{parametrixe}}

Let us recall that the goal of this proposition is to construct an operator $ \hat{P}_T $ on  $ (0,\infty) \times {\mathcal S} $ such that on one hand $ \hat{P} $ and $ \hat{P}_T $ coincide near infinity and on the other hand $ \hat{P}_T $ satisfy appropriate resolvent estimates ((\ref{1ertheoreme}), (\ref{rescalezero}) and (\ref{pourcompacite})).

\medskip

{\noindent \bf Construction of $ P_T $.} By Proposition \ref{reductiongeometriqueprop}, we assume that the metrics $ g_0 $ introduced in (\ref{exactconicalmetric}) and $ g (= a (r)dr^2 + 2r b (r)dr + r^2 h(r)  $, see (\ref{formedeG})) satisfy $ \mbox{div}_{g_0} = \mbox{div}_g $ near infinity.
More precisely, for any vector field $ V $ on $ (R,\infty) \times {\mathcal S} $, which we split as $ V = (V_1,V^{\mathcal S}) $ using the isomorphism $ T ((R,\infty) \times {\mathcal S} ) \approx 
T (R,\infty) \times T {\mathcal S} $, we have
$$ \mbox{div}_{g} (V) = \mbox{div}_{g_0} (V) = \frac{1}{r^{n-1}} \frac{\partial}{\partial r} \big( r^{n-1} V_1 \big) + \mbox{div}_{h_0} (V^{\mathcal S}) . $$
We then recall that  $ - \Delta_{g} $ (that is $\kappa^* \hat{P}  \kappa_*$ near infinity) is for $ r \gg 1 $ of the form
$$ -\Delta_{g} u = -\Delta_{g_0} u - \mbox{div}_{g_0} \big( K_G^{\rm sc} u\big)  , $$
for some tensor $ K_G $ such that,  by (\ref{hypothesedecroissance}),
\begin{eqnarray}
 || \partial_r^k K_G ||_{L^{\infty}} \leq C_k r^{-\rho-k} , \qquad k \geq 0 . \label{moinsrho}
\end{eqnarray}
(See also after (\ref{definitionPK}) for the notation sc.) We introduce $ \varphi_0 \in C^{\infty} (\Ra) $ such that $ \varphi_0 \equiv 1 $ on $ [1,\infty ) $ and $ \mbox{supp}(\varphi_0) \subset [1/2 , \infty ) $. Then, we define $$ T = \varphi_0 (r/R) K_G , $$
which is $ N $-admissible for all $ N $ by (\ref{moinsrho}), and set
\begin{eqnarray}
 P_T  u: = -\Delta_{g_0} u - \mbox{div}_{g_0} \big( T^{\rm sc} u \big) . \label{pourestimationpedagogie}
\end{eqnarray} 
We let $ \hat{P}_T $ be the associated self-adjoint realization defined according to Proposition \ref{pourfirstitem}. So defined, $ \hat{P}_T $ satisfies the item 1 of Proposition \ref{parametrixe}. Furthermore, by (\ref{moinsrho}), we have $ ||| T |||_0 \lesssim R^{-\rho} $ and $ ||| T |||_1 \lesssim R^{-\rho} $ (see Definition \ref{aumoinspourlanorme} for the norms $ ||| \cdot |||_k $), hence these norms are as small as we wish by choosing $R$ large enough. This  allows us to use the results of Propositions \ref{Mourrenonblackbox} and \ref{Ptrav} to prove the estimates (\ref{1ertheoreme}), (\ref{rescalezero}) and (\ref{pourcompacite}) as follows.

\bigskip

\noindent {\bf Proof of (\ref{1ertheoreme})} Let $ z = \lambda + i \delta $ with $ \delta \in \Ra \setminus \{ 0 \} $. Assume first that $ \lambda > 0 $. Setting $ \delta^{\prime} = \delta / \lambda $, it is straightforward  that
\begin{eqnarray}
 r^{-1} (\hat{P}_T - z)^{-1} r^{-1} & = & (\lambda^{1/2} r)^{-1} (\lambda^{-1} \hat{P}_T - 1 - i \delta^{\prime})^{-1} (\lambda^{1/2} r)^{-1} \nonumber \\
 & = & e^{itA}  r^{-1} ( \hat{P}_{T^t} - 1 - i \delta^{\prime})^{-1}  r^{-1} e^{-itA} , \label{rescalingfinal} 
\end{eqnarray} 
by choosing $ t = \ln (\lambda^{1/2}) $ (see (\ref{dilations}) for $ e^{itA} $). Notice that $ ||| T |||_k = ||| T^t |||_k $ for all $k$.  We next write
\begin{eqnarray}
 r^{-1} & = & r^{-1} (\alpha A + i) (\alpha A + i)^{-1} \nonumber \\
 & = & B (\alpha A + i)^{-1}, \label{sousexploite}
\end{eqnarray}
where $ B := \alpha D_r +  (\alpha n / 2 i + i) r^{-1} $ is bounded from $ H_0^1 $ to $ L^2 $ by the Hardy inequality (\ref{Hardy1}). This implies that
\begin{eqnarray*}
|| r^{-1} (\hat{P}_T - z)^{-1} r^{-1} ||_{L^2 \rightarrow L^2} 
 & = & || B (\alpha A+i)^{-1} ( \hat{P}_{T^t} - 1 - i \delta^{\prime})^{-1} (\alpha A-i)^{-1} B^* ||_{L^2 \rightarrow L^2} , \\
 & \leq & || B ||_{H_0^1 \rightarrow L^2}^2 || (\alpha A+i)^{-1} ( \hat{P}_{T^t} - 1 - i \delta^{\prime})^{-1} (\alpha A-i)^{-1}  ||_{H^{-1} \rightarrow H^1}
\end{eqnarray*} 
so the result follows from Proposition \ref{Mourrenonblackbox}. When $ \lambda \leq 0 $, the proof is even simpler and does not use Proposition (\ref{Mourrenonblackbox}). It suffices to use the Hardy inequality (\ref{Hardy1}) to see that
$$ || r^{-1} (\hat{P}_K - z)^{-1} r^{-1} ||_{L^2 \rightarrow L^2} \leq C || \hat{P}_K^{1/2} (\hat{P}_K - z)^{-1} \hat{P}_K^{1/2} ||_{L^2 \rightarrow L^2} , $$
whose right-hand side is bounded uniformly in $z$ by the spectral theorem.
\finpreuve

\bigskip

To prove (\ref{rescalezero}) and (\ref{pourcompacite}), we still need two technical lemmas.

\begin{lemm} \label{pourderivee} Fix $ s \geq 0 $  and $ 0 < \alpha < 1 $. Then there exist $ \varrho > 0 $ and $ C > 0 $ such that
$$ || \scal{r}^{-s-2} (\hat{P}_{K}+1)^{-1} u ||_{L^2 ({\mathcal M}_0)} \leq C || \scal{r}^{-s-1} (\alpha A+i)^{-1} u ||_{L^2 ({\mathcal M}_0)} , $$
for all $ u \in L^2 $ and all $K$ such that $ ||| K |||_{0} \leq \varrho $.
\end{lemm}

\noindent {\it Proof.} Without loss of generality we can replace $ \scal{r}$ by $ \scal{\delta r} $, with $ \delta > 0 $ small enough to be fixed below. 
We next remark that $ \bar{P}_0 + \bar{I} $ is an isomorphism from $ H_0^1 $ to $ H^{-1} $ (recall that $ \bar{P}_0 $ is the $ H_0^1 \rightarrow H^{-1} $ closure of $ - \Delta_{g_0} $) since $( (\bar{P}_0 + \bar{I}) u , v ) =  ( u,v)_{H_0^1} $.
%
%
Therefore, by Proposition \ref{propositionalgebreanalyse},  $ \bar{P}_K  + \bar{I} $ is also such an isomorphism (with norm in a fixed neighborhood of $ 1 $) if $ ||| K |||_0 $ is small enough. 
If $ \nu \geq 0 $ is a fixed real number, the operator 
\begin{eqnarray}
P_{K,\delta,\nu} = \scal{\delta r}^{-\nu} P_{K} \scal{\delta r}^{\nu} , \label{rappel1}
\end{eqnarray} 
defined first on $ C_0^{\infty} $, has a bounded closure $ \bar{P}_{K,\delta,\nu} $ to $ H_0^1 $ 
such that
$$ ||  \bar{P}_{K,\delta,\nu} -  \bar{P}_{K} ||_{H_0^1 \rightarrow H^{-1}} \leq C (1+||| K |||_0)  \delta , $$
where the constant is independent of $ \nu $ as long as $ \nu $ belongs to a bounded set.
This is easily seen from (\ref{definitionformeK}) and (\ref{definitionPK}), the factor $ \delta $ coming from commutations between $ \scal{\delta r}^{\pm \nu} $ and $ \partial_r $. For $ \delta $ small enough, 
 $ \bar{P}_{K,\delta,\nu} + \bar{I} $ is also an isomorphism between $ H_0^1 $ and $ H^{-1} $ and by construction (plus a routine verification which we omit), we obtain
\begin{eqnarray}
 \scal{ \delta r}^{-\nu} (\bar{P}_K + \bar{I})^{-1} = (\bar{P}_{K,\delta,\nu} + \bar{I})^{-1} \scal{\delta r}^{-\nu} , \label{decroissancepreservee}
\end{eqnarray}
as operators from $ H^{-1} $ to $ H_0^1 $. Here $ \scal{\delta r}^{-\nu} $ acts on $ H^{-1} $ in the distributions sense\footnote{{\it i.e.} $ ( \scal{\delta r}^{-\nu} E , u ) = (E, \scal{\delta r}^{-\nu} u) $ for all $ E \in H^{-1} $ and $ u \in H_0^1 $}. By composition with $ \bar{I} $ and  the item 2 of Proposition \ref{proplongue}, we get
\begin{eqnarray}
  \scal{ \delta r}^{-\nu} (\hat{P}_K + 1)^{-1}  =  (\bar{P}_{K,\delta,\nu} + \bar{I})^{-1} \bar{I} \scal{\delta r}^{-\nu} .  \label{rappel2}
\end{eqnarray}  
Using this identity with $ \nu = s + 2 $, we can thus write
$$ \scal{ \delta r}^{-s-2} (\hat{P}_K + 1)^{-1} = (\bar{P}_{K,\delta,\nu} + \bar{I})^{-1} \bar{I} \scal{\delta r}^{-1} (\alpha A + i) (\alpha A+i)^{-1} \scal{\delta r}^{-s-1} , $$
where, using that $ \scal{\delta r}^{-1} (\alpha A + i) $ maps $ L^2 $ in $ H^{-1} $, we thus obtain
$$ || \scal{ \delta r}^{-s-2} (\hat{P}_K + 1)^{-1} u ||_{L^2} \leq C || (\alpha A+i)^{-1} \scal{\delta r}^{-s-1} u ||_{L^2} . $$
  To swap the positions of $ (\alpha A+i)^{-1} $ and $ \scal{\delta r}^{-s-1} $, we write
 $$ (\alpha A+i)^{-1} \scal{\delta r}^{-s-1} = (\alpha A+i)^{-1} \scal{\delta r}^{-s-1} (\alpha A + i) (\alpha A + i)^{-1} , $$
 and observe that $ (\alpha A+i)^{-1} \scal{\delta r}^{-s-1} (\alpha A + i) = B (\alpha ,s,\delta) \scal{\delta r}^{-s-1} $ for some bounded (and explicitly computable) operator $ B (\alpha,s,\delta) $. The result follows.  \finpreuve

\bigskip

\begin{lemm} \label{pourrescalee} Fix $ M , N \geq 0 $. There exist $ \alpha_0 > 0 $, $ \varrho > 0 $ and $ C > 0 $ such that, for all integer $ 0 \leq k \leq N-1 $,
$$ || r^{-1} \scal{r}^{-k} (\hat{P}_K +1)^{-1-N} u ||_{L^2 ({\mathcal M}_0)} \leq C ||  (\alpha A+i)^{-1-k} u ||_{L^2 ({\mathcal M}_0)} , $$
for all $ u \in L^2 $, all $ 0 < \alpha < \alpha_0 $ and all $K$ such that $ ||| K |||_{0} \leq \varrho $, $ |||K|||_N \leq M $.
\end{lemm}

We state this result for general admissible perturbations $ K $ but we will apply it with $ K = T^t $,  using that  $ || T^t ||_k = || T ||_k $ for all integer $k $ and all  $ t \in \Ra $.

\bigskip

\noindent {\it Proof.}  We will use (\ref{rappel1}) and (\ref{rappel2}) from the proof of Lemma \ref{pourderivee}. 
By iteration of (\ref{rappel2}), we obtain on one hand
\begin{eqnarray}
\scal{\delta r}^{-k} (\hat{P}_K+1)^{-k-1} & = & \big( \bar{P}_{K,\delta,k} + \bar{I} \big)^{-1} \bar{I} \scal{\delta r}^{-k} (\hat{P}_K+1)^{-k} \nonumber \\
 & = &  \big( \bar{P}_{K,\delta,k} + \bar{I} \big)^{-1} \bar{I} \scal{\delta r}^{-1} \big(\bar{P}_{K,\delta,k-1}+1 \big)^{-1} \bar{I} \scal{\delta r}^{-(k-1)} (\hat{P}_K+1)^{-(k-2)} \nonumber \\
 & = &  \big( \bar{P}_{K,\delta,k} + \bar{I} \big)^{-1} \bar{I} \prod_{\nu=k-1}^{0} \scal{\delta r}^{-1} \big(\bar{P}_{K,\delta,\nu}+1 \big)^{-1} \bar{I} , \label{pedagogie}
\end{eqnarray}
where the product is the composition from the left to the right decreasingly in $\nu$.
 On the other hand, for any $ \nu \geq 0 $ and an integer $ j \geq 0 $, we can consider the operator
\begin{eqnarray}
 P_{K,\delta,\nu}^{j,\alpha}  := (\alpha A + i)^{-j}P_{K,\delta,\nu}  (\alpha A + i)^j , \label{theoperator}
\end{eqnarray}
on $ C_0^{\infty} $. This can be written as the sum of $ P_{K,\delta,\nu}  $ and a linear combination of nonnegative powers of $ (\alpha A + i)^ {-1} $ composed with commutators of $ P_{K,\delta,\nu} $ and $ \alpha A $. It follows that
$$ \big| \big| ( P_{K,\delta,\nu}^{j,\alpha} - P_{K,\delta,\nu} ) u \big| \big|_{H^{-1}} \leq C \alpha || u ||_{H_0^1} , \qquad u \in C_0^{\infty}, 0 < \alpha \ll 1 . $$ 
This implies that $ P_{K,\delta,\nu}^{j,\alpha} $ has a closure $ \bar{P}_{K,\delta,\nu}^{j,\alpha} $ to $ H_0^1 $. If $ \alpha $ is small enough, $  \bar{P}_{K,\delta,\nu}^{j,\alpha} + \bar{I} $ is an isomorphism between $ H_0^1 $ and $ H^{-1} $ since $ \bar{P}_{K,\delta,\nu} $ is  for $ \delta $ small enough (cf the proof of Lemma \ref{pourderivee}). Moreover, (\ref{theoperator}) implies 
\begin{eqnarray}
 (\alpha A+i)^{-j} \big( \bar{P}_{K,\delta,\nu} + \bar{I} \big)^{-1} = \big( \bar{P}_{K,\delta,\nu}^{j,\alpha} + \bar{I} \big)^{-1} (\alpha A+i)^{-j}  .
 \label{wecanuse} 
\end{eqnarray} 
This is formally obvious but requires an argument since we cannot obviously compose both sides of (\ref{theoperator}) with $ (\alpha A + i)^{-j} $ for this does not preserve $ C_0^{\infty} $ in general. To justify this formula we use the Lemma 3.1 of \cite{BoucletCPDE} as in the proof of Proposition \ref{Ptrav}, namely that for any $ v \in H_0^1 $ we can find a sequence $ (u_m)_{m \in \Na} $ of $ C_0^{\infty} $ such that
$$ (\alpha A + i)^j u_m \rightarrow v \qquad \mbox{and} \qquad u_m \rightarrow (\alpha A + i)^{-j} v $$ 
in $ H_0^{1} $. Then (\ref{theoperator}) yields $ (\alpha A + i)^{-j} \bar{P}_{K,\delta,\nu} v = \bar{P}_{K,\delta,\nu}^{j,\alpha} (\alpha A + i)^{-j} v $ which then implies (\ref{wecanuse}). The interest of (\ref{pedagogie}) and (\ref{wecanuse}) is the following one.
%
%
%
%
After multiplication by $ r^{-1} $, we  write in the last line of (\ref{pedagogie})
\begin{eqnarray*}
r^{-1} \big( \bar{P}_{K,\delta,k} + \bar{I} \big)^{-1} \bar{I} & = & r^{-1} (\alpha A+i) (\alpha A+i)^{-1} \big( \bar{P}_{K,\delta,k} + \bar{I} \big)^{-1} \bar{I} \\
& = & r^{-1} (\alpha A + i) \big( \bar{P}_{K,\delta,k}^{1,\alpha } + \bar{I} \big)^{-1} \bar{I} (\alpha A+i)^{-1} .
\end{eqnarray*} 
The operator $ (\alpha A+i)^{-1} $ in the right-hand side, which falls on the operator $ \scal{\delta r}^{-1} \big(\bar{P}_{K,\delta,k-1}+1 \big)^{-1} $, is then rewritten as $ (\alpha A+i) (\alpha A + i )^{-2} $ so that we can use  (\ref{wecanuse}) with $ j = 2 $ and $\nu = k-1 $. By iteration,
 we see that $ r^{-1} \scal{\delta r}^{-k} (\hat{P}_K+1)^{-k-1} $ reads 
 $$ r^{-1} (\alpha A+i) \big( \bar{P}_{K,\delta,k}^{1,\alpha } + \bar{I} \big)^{-1} \bar{I} \left( \prod_{l=k-1}^{0} B(k,l) (\alpha A+i) \scal{\delta r}^{-1} \big(\bar{P}_{K,\delta,l}^{k+1-l,\alpha}+1 \big)^{-1} \bar{I}   \right) (\alpha A+i)^{-k-1} , $$
 with $ B(k,l) $   such that $ (\alpha A+i)^{-(k+1-l)} \scal{\delta r}^{-1} (\alpha A+i)^{k+1-l} = B (k,l) \scal{\delta r}^{-1} $. Each  $ B(k,l) $ is clearly bounded on $ L^2 $.
 By using Lemma \ref{injectiondomaine} and the fact that $ r^{-1} (\alpha A+i) $ is bounded on $ H_0^1 $, we conclude that
\begin{eqnarray*}
 || r^{-1} \scal{\delta r}^{-k} (\hat{P}_K+1)^{-N-1} u ||_{L^2} & \leq & C || (\alpha A+i)^{-k-1} (\hat{P}_K + 1)^{-(N-k)} u ||_{L^2} \\
 & \leq & C || (\alpha A + i)^{-k-1} u ||_{L^2}
\end{eqnarray*} 
where, to get the second line, we used (\ref{wecanuse}) with $ j = k+1 $ and $ \nu = 0 $ (in this case we have $ (\hat{P}_K + 1)^{-1} = (\bar{P}_{K,0,0} + \bar{I} )^{-1} \bar{I} $).  This completes the proof.  \finpreuve


\bigskip

\noindent {\bf Proof of (\ref{rescalezero})} We start with a general remark. By iterating the resolvent identity for any self-adjoint operator $H \geq 0$, we have
\begin{eqnarray}
 (H - \zeta)^{-1} = \sum_{j=0}^{2N}  (1+\zeta)^j (H+1)^{-1-j}  + (\zeta+1)^{2N+1} (H +1)^{-N}(H - \zeta )^{-1} (H+1)^{-N} , \label{localisationfaible}
\end{eqnarray} 
for all $ \zeta \in \Ca \setminus \Ra $. 
By differentiating $k$ times in $\zeta$, we see that if $\zeta$ belongs to a bounded subset of $ \Ca \setminus \Ra $, there exists $ C > 0 $ such that, for all bounded operator $ W $ with operator norm at most $1$, 
\begin{eqnarray}
 || W (H-\zeta)^{-1-k} W || \leq C \left( 1 + \sum_{l \leq k} || W (H+1)^{-N} (H-\zeta)^{-1-l} (H+1)^{-N} || \right) . \label{identiteresolventeMourre}
\end{eqnarray} 
Here $ || \cdot || $ is the operator on the Hilbert space where $ H $ is defined. From now on, we consider $ \hat{P}_T $  and $|| \cdot || $ is the operator norm on $ L^2 ({\mathcal M}_0)$. We let $ \epsilon^2 Z = \lambda (1 + i \delta^{\prime}) $, with $ \delta^{\prime} \in \Ra \setminus \{ 0 \} $ (note that $ \lambda \sim \epsilon^2 $). Then
\begin{eqnarray}
 (\epsilon^{-2} \hat{P}_T - Z)^{-1-k} = \mbox{Re}(Z)^{-1-k} \big( \lambda^{-1} \hat{P}_T - 1 - i \delta^{\prime} \big)^{-1-k} . 
\label{algebrevisible}
\end{eqnarray} 
 Similarly to (\ref{rescalingfinal}), we consider the family of rescaled operators $ \hat{P}_{T^t} = e^{itA} ( \lambda^{-1} \hat{P}_T )  e^{-itA} $. We observe that
 \begin{eqnarray*}
  || (\epsilon r)^{-1} \scal{\epsilon r }^{-k} \big( \epsilon^{-2} \hat{P}_T - Z  \big)^{-1-k} \scal{\epsilon r }^{-k}( \epsilon r)^{-1} ||
  \lesssim 
  ||  r^{-1} \scal{ r }^{-k} \big( \hat{P}_{T^t} - 1 - i \delta^{\prime} \big)^{-1-k} \scal{ r }^{-k} r^{-1} || .
\end{eqnarray*} 
Indeed the right-hand side equals $ || (\lambda^{1/2} r)^{-1} \scal{\lambda^{1/2} r }^{-k} \big( \lambda^{-1} \hat{P}_T - 1 - i \delta^{\prime} \big)^{-1-k} \scal{\lambda^{1/2} r }^{-k}(\lambda^{1/2} r)^{-1}|| $ by rescaling, and  this quantity is bounded from above and below by the left-hand side, using (\ref{algebrevisible}) and the fact that $ \lambda / \epsilon^2 $ belongs to a compact see of $ (0,\infty) $. Then, by using (\ref{identiteresolventeMourre}) and Lemma \ref{pourrescalee}, we have
$$  ||  r^{-1} \scal{ r }^{-k} \big( \hat{P}_{T^t} - 1 - i \delta^{\prime} \big)^{-1-k} \scal{ r }^{-k} r^{-1} || \lesssim 1 +
\sum_{l \leq k} || (\alpha A+i)^{-1-l} (\hat{P}_{T^t}-1-i\delta^{\prime})^{-1-l} (\alpha A-i)^{-1-l} ||  $$
so the result follows from Proposition \ref{Mourrenonblackbox}. \finpreuve


\bigskip

\noindent {\bf Proof of (\ref{pourcompacite})}  By using (\ref{identiteresolventeMourre}) with $z = \zeta$, we obtain (as long as $z$ is bounded)
$$ || \scal{r}^{-2-s} (\hat{P}_T-z)^{-2} \scal{r}^{-2-s} || \lesssim 1 + \sum_{k=1}^2 || \scal{r}^{-2-s} (\hat{P}_T + 1)^{-1} (\hat{P}_T-z)^{-k} (\hat{P}_T + 1)^{-1} \scal{r}^{-2-s} || . $$
The term corresponding to $ k = 1 $ is clearly bounded for $ |\mbox{Re}(z)| \leq 1 $ and $0 < |\mbox{Im}(z) | \leq 1 $ by using (\ref{1ertheoreme}) and the fact that $ (\hat{P}_T+1)^{-1} $ preserves the decay $ \scal{r}^{-2-s} $ (see (\ref{decroissancepreservee})). Therefore, it suffices to consider the term corresponding to $ k =2 $. By Lemma \ref{pourderivee}, this term is controlled by
$$ 
 || (r+1)^{-1-s} (\alpha A+i)^{-1} (\hat{P}_T-z)^{-2} (\alpha A-i)^{-1} (r+1)^{-1-s} ||_{L^2 \rightarrow L^2} 
. $$
We assume first that $ \mbox{Re}(z)=: \lambda  $ is positive.
By using the same rescaling as in (\ref{rescalingfinal}), the above norm reads
\begin{eqnarray}
   \mbox{Re}(z)^{s-1} || (r + \lambda^{1/2})^{-1-s} (\alpha A+i)^{-1} (\hat{P}_{T^t}-1-i\delta^{\prime})^{-2} (\alpha A-i)^{-1} (r+\lambda^{1/2})^{-1-s} ||_{L^2 \rightarrow L^2}  . \label{abovenorm} 
\end{eqnarray}
  Therefore, it suffices to show that the norm in (\ref{abovenorm}) is bounded uniformly in $ \delta^{\prime} $ and $ \lambda $ (recall that $t = \ln( \lambda^{1/2}) $). Using (\ref{sousexploite}), we write 
\begin{eqnarray*}
 ( r + \lambda^{1/2})^{-1-s} (\alpha A+i)^{-1} & = &  ( r + \lambda^{1/2} )^{-s}  \left( \frac{\alpha r}{r + \lambda^{1/2}} D_r + \frac{\alpha n / 2i + i}{r + \lambda^{1/2}} \right) (\alpha A+i)^{-2}  \\
 & = : & ( r + \lambda^{1/2} )^{-s}  B_{\lambda} (\alpha A+i)^{-2}
\end{eqnarray*}
so  that the norm in (\ref{abovenorm})  equals precisely 
$$ \big| \big| ( r+ \lambda^{1/2})^{-s} B_{\lambda} (\alpha A+i)^{-2} (\hat{P}_{T^t}-1-i\delta^{\prime})^{-2} (\alpha A-i)^{-2} B^*_{\lambda} ( r + \lambda^{1/2} )^{-s} \big| \big|_{L^2 \rightarrow L^2} . $$
By Lemma \ref{pour2eme} and Proposition \ref{Ptrav}, this norm is bounded by a constant (independent of $ \lambda $ and $ \delta^{\prime} $) times the product of the following powers of norms
\begin{eqnarray*}
\big| \big|  (\alpha A+i)^{-2} \hat{P}_{T^t}^2 (\hat{P}_{T^t}-1-i\delta^{\prime})^{-2} (\alpha A-i)^{-2}  \big| \big|_{L^2 \rightarrow L^2}^{1/4} \\
\big| \big|  ( r+\lambda^{1/2} )^{-2s} (\alpha A+i)^{-2}  (\hat{P}_{T^t}-1 \pm i  \delta^{\prime})^{-2} \hat{P}_{T^t} (\alpha A-i)^{-2}  \big| \big|_{L^2 \rightarrow L^2}^{1/4} \\
\big| \big| ( r + \lambda^{1/2} )^{-2s} (\alpha A+i)^{-2}  (\hat{P}_{T^t}-1-i\delta^{\prime})^{-2} (\alpha A-i)^{-2} (r + \lambda^{1/2})^{-2s} \big| \big|_{L^2 \rightarrow L^2}^{1/4} .
\end{eqnarray*}
Since $ s < 1/2 $, the Hardy inequality allows to drop the weight $ ( r + \lambda^{1/2} )^{-2s} $  in the second and third lines up to the replacement of the $ L^2 \rightarrow L^2 $ norm
by the $ H^{-1} \rightarrow H_0^{1} $ norm. The uniform boundedness of these norms then follows from Proposition \ref{Mourrenonblackbox}. The proof in the case $ \mbox{Re}(z)< 0 $ is similar; one only has to replace $ (\hat{P}_{T^t} - 1 - i  \delta^{\prime}) $ by $(\hat{P}_{T^t} + 1 - i  \delta^{\prime})$ so that we do not need to use Proposition \ref{Mourrenonblackbox}. This completes the proof. \finpreuve


\end{document}